\documentclass{article}
 
\usepackage{PRIMEarxiv}

\usepackage{microtype}
\usepackage{graphicx}
\usepackage{subfigure}
\usepackage{booktabs} 
\usepackage{siunitx}

\usepackage[
	colorlinks=true,
	hypertexnames=false
	]{hyperref}

\hypersetup{
  colorlinks   = true, 
  urlcolor     = magenta, 
  linkcolor    = blue, 
  citecolor   = teal 
}

\usepackage{comment}

\usepackage{pgfplots}
\DeclareUnicodeCharacter{2212}{-}
\usepgfplotslibrary{groupplots,dateplot}
\usetikzlibrary{patterns,shapes.arrows}
\pgfplotsset{compat=newest}
\pgfplotsset{every tick label/.append style={font=\small}}
\usepackage{caption}

\usepackage{amsmath,amsfonts,amssymb,amsthm,bbm,graphicx,enumerate,times,mathdots,braket}
\usepackage[capitalise,compress]{cleveref}
\usepackage{mathtools}
\usepackage{tikz}
\usetikzlibrary{positioning}
\usepackage{graphicx}
\usepackage[normalem]{ulem}
\usepackage{float}
\usepackage{comment}
\usepackage{csquotes}
\usepackage{url}
\usepackage{enumerate}  
\usepackage{cuted}

\usepackage[ruled]{algorithm2e}
\usetikzlibrary{decorations.pathreplacing}
\newtheorem{theorem}{Theorem}

\newtheorem{definition}[theorem]{Definition}

\newtheorem{corollary}[theorem]{Corollary}
\DeclareMathOperator{\sign}{sign}
\DeclareMathOperator{\id}{\mathrm{Id}}
\DeclareMathOperator{\tr}{tr}
\DeclareMathOperator{\vspan}{span}
\DeclareMathOperator{\diag}{diag}
\DeclareMathOperator{\image}{Image}
\DeclareMathOperator{\rank}{rank}
\DeclareMathOperator{\poly}{poly}

\definecolor{jonas}{rgb}{.4,0,.7}
\definecolor{ingo}{rgb}{.7,0,.3}

\title{\centering A quantum inspired approach to learning dynamical laws from data---block-sparsity and gauge-mediated weight sharing}
\author{
  J. Fuksa${^{\ast, 2}}$, M. G{\"o}tte${^{\ast, 1}}$, I. Roth${^3}$, and J. Eisert${^{2,4}}$ \vspace{5pt}\\
  1 \tu \\
  2 \fu \\
  3 Quantum Research Centre, Technology Innovation Institute,
  Abu Dhabi \\
  4 Fraunhofer Heinrich Hertz Institute, Germany \vspace{5pt}\\
  $\ast$ \textit{J.F. and M.G. have contributed equally and correspond under}\\
  \texttt{jonas.fuksa@fu-berlin.de} \\
}
\newcommand{\tu}{Institute of Mathematics, Technische Universit{\"a}t Berlin, Germany}
\newcommand{\fu}{Dahlem Center for Complex Quantum Systems, Freie Universit{\"a}t Berlin, Germany}

\begin{document}

\maketitle

\begin{abstract}
Recent years have witnessed an increased interest in recovering dynamical laws of complex systems in a largely data-driven fashion under meaningful hypotheses.
In this work, we propose a  scalable and numerically robust  method for this task, utilizing efficient block-sparse tensor train representations of dynamical laws, inspired by similar approaches in quantum many-body systems. Low-rank tensor train representations have been previously derived for dynamical laws of one-dimensional systems.
We extend this result to efficient representations of systems with $K$-mode interactions and  controlled approximations of systems with decaying interactions. 
We further argue that natural structure assumptions on dynamical laws, such as bounded polynomial degrees, can be exploited in the form of block-sparse support patterns of tensor-train cores. 
Additional structural similarities between interactions of certain modes can be accounted for by weight sharing within the ansatz. 
To make use of these structure assumptions, we propose a novel optimization algorithm, block-sparsity restricted alternating least squares with gauge-mediated weight sharing.
The algorithm is inspired by similar notions in machine learning and achieves a significant improvement in performance over previous approaches.
We demonstrate the performance of the method numerically on three one-dimensional systems -- the Fermi-Pasta-Ulam-Tsingou system, rotating magnetic dipoles and point particles interacting via modified Lennard-Jones potentials, 
observing a highly accurate and noise-robust recovery.  
\end{abstract}


\twocolumn

\section{Introduction}

Discovering dynamical laws that govern the time evolution of dynamical systems has been a central task in physics and engineering for centuries, from practical as well as fundamental perspective.
Historically, this task has been approached from two directions -- firstly by using expert physical knowledge and intuition, and secondly by using data obtained by measuring time evolution of the dynamical system in question. 
With the ever growing availability of large amounts of computational power and data, the second approach is becoming increasingly accessible \cite{schmidt_distilling_2009,brunton_discovering_2016,gels_multidimensional_2019,goesmann_tensor_2020,PhysRevLett.124.010508,CINDy,Kaheman_2022,Cornelio2021AIDC}.
However, since a crucial aspect of using data to learn dynamical laws is choosing the right hypothesis class, physical intuition about the system is essential in developing efficient data-driven algorithms.

One prominent recent approach is the \emph{sparse identification of non-linear dynamics} (SINDy) algorithm \cite{brunton_discovering_2016,PhysRevResearch.3.023255,PythonPackage}.
Here the learning task is phrased as a linear inversion problem for a chosen function dictionary.
To arrive at a physically motivated hypothesis class, the authors impose sparsity of the recovered dynamical laws with respect to this dictionary---implementing the principle of Occam's razor.

The remarkable success of the SINDy algorithm demonstrates that imposing structure in learning dynamical laws is immensely powerful. 
This is in spite of the fact that Occam's razor is a general heuristic principle that is not linked to any specific physical properties of the system at hand.

Building on the ideas of SINDy and its variant MANDy \cite{gels_multidimensional_2019}, in ref.~\cite{goesmann_tensor_2020} the authors proposed to use \emph{locality} in one dimension as the structure imposing physical principle.
They have shown that the resultant hypothesis class consists of low rank \emph{tensor trains} (TT) \cite{bachmayr_tensor_2016}, a specific type of an efficient tensor network representation of multivariate functions \cite{cohen_expressive_2015,stoudenmire_supervised_2016,levine_deep_2017}.
A similar result is widely known in the quantum many-body literature \cite{RevModPhys.93.045003,Orus-AnnPhys-2014,verstraete_renormalization_2004,haegeman_geometry_2014}, where it has been shown that low rank TTs (known in this context as \emph{matrix product states}) parametrize ground states of one-dimensional local Hamiltonians \cite{verstraete_matrix_2006}, as well as other states of physical importance \cite{eisert_colloquium_2010,schuch_entropy_2008,verstraete_matrix_2008}.

Another structural observation made by the quantum many-body community is that the TT representations of quantum states that are symmetric under the local action of some symmetry group have tensor cores with a characteristic sparse support, a property dubbed \emph{block-sparsity} \cite{singh_tensor_2010}.
Recently, it has been proven that functions with bounded polynomial degree (or an equivalent notion suitable for the chosen function dictionary) also admit block-sparse TT representations \cite{gotte_block-sparse_2021-1,bachmayr_particle_2021}.
Being able to control the polynomial degree of non-linear functions is in the context of learning dynamical laws a promising primitive.
On a technical level, limiting the total polynomial degree offers a natural truncation of a multivariate function space via Taylor's theorem.
On a conceptual level, high-degree polynomials in dynamical laws correspond to terms that vary quickly with respect to many modes.
The appearance of such terms in dynamical laws is generally a sign that a better choice of coordinates parametrizing the state space could be made.
Similar arguments also apply to other dictionaries, e.g., for trigonometric functions limiting the total ``degree'' corresponds to neglecting fast oscillating terms in the dynamical law.

Bounding a notion of total degree suitable for a chosen multi-variate function dictionary, block-sparse TTs arise as a natural, physically well-motivated efficient restriction of the ansatz class for dynamical laws. 
In another context, block-sparsity has been observed to improve the performance of tensor network optimization methods, in terms of computational, memory and sample complexity \cite{gotte_block-sparse_2021-1,bachmayr_particle_2021}.

The contribution of this work is three-fold:
(1) We broaden the range of physical principles that lead to efficient TT representations of dynamical laws.
In particular, we show that systems with $K$-mode interactions admit an efficient TT representation and that systems with algebraically decaying interactions can be approximated with bounded error by a system with efficient TT representation.
(2) We utilize block-sparsity in the context of learning dynamical laws---demonstrating a significantly improved performance of the resulting method compared to previous work.
(3) We use self-similarity between certain modes in the system to further restrict the search space and develop a new optimization algorithm for the resultant hypothesis class, referred to as \emph{ALS optimization with gauge mediated weight sharing (ALS-GMWS)}, which is inspired by similar notions in machine learning \cite{riemer_role_2020}.
We show numerically that (2) and (3) improve scalability of the method, allowing us to learn dynamical laws of systems more than three times larger than the systems presented in ref.~\cite{goesmann_tensor_2020}.

The remainder of this work is structured as follows: Our setting of learning dynamical laws is formally introduced in \cref{sec:setup}.
In \cref{sec:tensor networks} we give a brief primer into tensor networks.
\cref{sec:block sparsity} shows how a natural truncation of the function space leads to an efficient parametrization via block sparse TTs.
In \cref{sec:efficient} we prove that low rank TTs parametrize dynamical laws of systems obeying generic physical principles.
The concept of self-similarity and the ALS-GMWS algorithm is presented in \cref{sec:gauge mediated weight sharing}.
Finally, we numerically demonstrate the performance of our method in \cref{sec:numerics}, before concluding in \cref{sec:conclusions}.

\section{Setup}\label{sec:setup}

We will use the notation $[N] := \{1,\dots,N\}$ for any integer $N$.
Consider a dynamical system with state space $\mathcal{S} \subset \mathbb{R}^d$, such that its state is described by a real $d$-dimensional vector $x = (x^{(1)},\dots,x^{(d)})$. 
We call each of the degrees of freedom $x^{(i)}$ a \emph{mode}.
The time evolution is a smooth curve $x(t): \mathbb{R} \rightarrow \mathcal{S}$, which is, given initial conditions $x_0$, generated by
\begin{equation}\label{eq:dynamical law}
    \dot x(t) = f(x(t)) \quad \text{or} \quad \ddot x(t) = f(x(t)),
\end{equation}
where we have restricted our attention to time-independent systems.
We assume that it is known which of these forms governs the dynamical system of interest and we will refer to both eq.~\eqref{eq:dynamical law} and the function $f: \mathcal{S} \rightarrow \mathbb{R}^d$ as the \emph{dynamical law}.
The first form of eq.~\eqref{eq:dynamical law} can arise, e.g., from Hamilton's equations, where $f(x) = \{x,H\}$ with $H$ the Hamiltonian and $\{\cdot,\cdot\}$ the Poisson bracket, while the second form appears, e.g., in Newton's equations, where $f_k(x)$ is the total force acting on the $k$-th mode in the state $x$.
By \emph{learning a dynamical law} we mean identifying $f(x)$ given data pairs $(x_i,y_i) \in \mathbb{R}^d \times \mathbb{R}^d$ for $i \in [M]$ with the relation $y_i(x_i) \approx f(x_i)$. 
The data may be coming from time-series measurements of a trajectory with the gradients $y_i$ approximated, e.g., by the method of finite differences.

To learn $f$, we choose a \emph{function dictionary} $\{\Psi_i: \mathbb{R} \rightarrow \mathbb{R}\}_{i\in [p]}$ of linearly independent functions.
Forming a product basis 
\begin{equation}\label{eq:product basis}
    \Phi_{i_1,\dots,i_d}(x) = \Psi_{i_1}(x^{(1)})\dots \Psi_{i_d}(x^{(d)})
\end{equation}
we obtain the space 
\begin{equation}
    \mathcal{A} = \vspan\{\Phi_{i_1,\dots,i_d}: \, i_k \in [p] \, \forall k \in [d]\},
\end{equation}
with elements mapping $\mathbb{R}^d$ to $\mathbb{R}$, so that $\mathcal{A}^d$ can be used as a search space for $f: \mathbb{R}^d \rightarrow \mathbb{R}^d$.

Elements $g(x) \in \mathcal{A}$ are labeled by tensors $\Theta \in \mathbb{R}^{p^d}$ via 
\begin{equation}\label{eq:product basis decomp}
    g(x) = \sum_{i_1,\dots,i_d \in [p]^d} \Theta_{i_1,\dots,i_d} \Phi_{i_1,\dots,i_d}(x).
\end{equation}
Hence, the dimension of $\mathcal{A}$ is exponential in the system size $d$, limiting the scalability of learning a function within $\mathcal{A}$.
This is an example of the notorious \emph{curse of dimensionality}.
In order to overcome the curse, we use structural constraints on $\Theta$ due to known physical properties of the system to identify a physically relevant subspace $\mathcal{H} \subset \mathcal{A}$ with $\dim(\mathcal{\mathcal{H}}) \in O(\operatorname{poly}(d))$, defining the \emph{physical corner} of solutions.
Specifically, in this work, we show that for many dynamical systems of interest a suitable $\mathcal{H}$ is the set of low rank block-sparse TTs, as defined in \cref{sec:tensor networks}.

\section{A primer to tensor networks}\label{sec:tensor networks}

To set the stage and notation, we start with a brief introduction to tensor networks.
For a more thorough treatment of this topic, see ref.~\cite{bachmayr_tensor_2016,bridgeman_hand-waving_2017-1}.
The curse of dimensionality, as discussed in the previous section, renders high-order tensors hard to work with.
Not only are they hard to optimize over, even storing them and performing basic operations on them quickly becomes impractical as the order grows.
\emph{Tensor networks} are a way of decomposing high-order tensors in terms of contractions of smaller tensors.
Limiting the ranks of these contractions we define subspaces of tensor spaces with dimensions scaling polynomially in the order, averting the curse of dimensionality.
Crucially for us, such subspaces turn out to contain physically relevant tensors in many cases of interest.

Tensor networks can be visually represented by \emph{tensor network diagrams}.
Suppose an order $d$ tensor $T$, which is a contraction of tensors $\{C_i\}_i$, called \emph{tensor cores}.
To draw a tensor network diagram for $T$, we draw a node for each tensor core $C_i$.
For each of the indices of the tensor cores we draw an edge and connect those edges that represent indices that are contracted over in $T$.
The indices represented by connected edges are referred to as \emph{virtual indices}.
The remaining $d$ unconnected edges represent the $d$ \emph{physical indices} of $T$.

A widely studied example of such a network is the \emph{tensor train}, which is represented by the tensor network
\begin{equation}\label{eq:tensor_train}
\begin{tikzpicture}
    \node at (-1.4,0) {$T_{i_1,\ldots ,i_d} =$};
    \node [above] at (0,0) {$C_{1}$};
    \draw (0,-0.5) -- (0,0) -- (1,0);
    \draw [fill] (0,0) circle (0.5ex);
    \node [below] at (0,-0.5) {$i_1$};
    
    \node [above] at (1,0) {$C_{2}$};
    \draw (1,-0.5) -- (1,0) -- (1.5,0);
    \draw [fill] (1,0) circle (0.5ex);
    \node [below] at (1,-0.5) {$i_2$};

    \node at (2,0) {$\dots$};
    
    \node [above] at (3,0) {$C_{d-1}$};
    \draw (2.5,0) -- (4,0) -- (4,-0.5);
    \draw (3,0) -- (3,-0.5);
    \draw [fill] (3,0) circle (0.5ex);
    \node [below] at (3,-0.5) {$i_{d-1}$};

    \node [above] at (4,0) {$C_{d}$};
    \draw [fill] (4,0) circle (0.5ex);
    \node [below] at (4,-0.5) {$i_d$};

    \node [below] at (4.5,0) {.};
\end{tikzpicture}
\end{equation}
This contraction of tensors can be written in matrix notation as
\begin{equation}
    \begin{split}
    T_{i_1,\dots,i_d} =
    \mkern-27mu\sum_{\substack{j_1,\dots,j_{d-1} \\\in [r_1] \times \ldots \times [r_{d-1}]}} 
    \mkern-27mu(C_1)_{j_1}^{i_1} (C_2)_{j_1,j_2}^{i_2} \dots 
    (C_d)_{j_{d-1}}^{i_d},
    \end{split}
\end{equation}
where the superscripts on the tensor cores are the physical indices, while the subscripts are the virtual indices.
Any tensor $T \in \mathbb{R}^{d^p}$ can be represented as a TT, if we allow $r_k$ to scale exponentially with $d$. 
To obtain an efficient ansatz class that does not suffer from the curse of dimensionality, we have to bound $r_k$ by a polynomial in $d$.
We call $r = \min_{k\in[d]}r_k$ the TT rank of such a decomposition of $T$.

An important observation is that the TT decomposition is not unique.
The \emph{gauge transformation}
\begin{equation}\label{eq:gaugetrafo}
    \begin{split}
        C_1 &\mapsto C_1 A_1^{-1}, \qquad C_d \mapsto A_{d-1} C_d,\\
        C_\ell &\mapsto A_{\ell-1} C_\ell A_\ell^{-1} \quad \forall \ell \in \{2,\dots,d-1\}\,,        
    \end{split}
\end{equation}
applying a matrix multiplication with invertible matrices $A_\ell$ on the virtual indices of the tensor cores, leaves the TT invariant.
Part of this \emph{gauge freedom} is removed by imposing the so-called \emph{left canonical condition}
\begin{equation}\label{eq:left canonical}
    \begin{tikzpicture}[
        dot/.style = {circle,fill,minimum size=#1,inner sep=0pt,outer sep=0pt},
        dot/.default = 5pt
        ]
        \node (origin) at (0,0) {};
        \node[dot,label=above:$C_i$] [above=.15cm of origin] (c1) {};
        \node[dot,label=below:$C_i^*$] [below=.5cm of c1] (c1*) {};
        \node [left=.2 of c1] (dl) {};
        \node [right=.3 of c1] (dr) {};
        \node [left=.2 of c1*] (dl*) {};
        \node [right=.3 of c1*] (dr*) {};
        \node [right=.5 of origin] (eq) {$=$};
        \node [right=.5 of eq] (centerr) {};

        \draw (c1) -- (c1*);
        \draw (dr) -- (dl.center) -- (dl*.center) -- (dr*);

        \node (center) [right=1.7 of origin] {};
        \node [dot,label=left:$C_1$] (c) [above=.15 of center] {};
        \node [dot,label=left:$C_1^*$] (c*) [below=.5 of c] {};
        \node (cr) [right=.3 of c] {};
        \node (c*r) [right=.3 of c*] {};
        \node (eq2) [right=.5 of center] {$=$};

        \draw (c*r) -- (c*) -- (c) -- (cr);

        \node (centerr2) [right=.0 of eq2] {};
        \node (top2) [above=.15 of centerr2] {};
        \node (bottom2) [below=.5 of top2] {};
        \node (rightt) [right=.3 of top2] {};
        \node (rightb) [right=.3 of bottom2] {};

        \draw (rightb) -- (bottom2.center) -- (top2.center) -- (rightt);
    \end{tikzpicture}
\end{equation}
$\forall i \in \{2,\dots,d-1\}$.
The remaining gauge freedom can be shown to be a unitary transformation.

\section{Block-sparsity}\label{sec:block sparsity}

In choosing the finite local dictionary, we are truncating the univariate function space, e.g. bounding the polynomial degree or a similar notion relevant to the dictionary at hand.
However, if we now take the $d$-fold tensor product $\mathcal{A}$ of the univariate dictionaries as the search space for a multivariate $f_k$, we introduce  terms that are of higher total degree than the imposed local truncation. 
Thus, from the perspective of multivariate Taylor series, the multivariate dictionary appears inconsistently truncated. 
Here, it is more natural to work with complete function spaces of bounded \emph{total degree}.
In this section we show that such a natural restriction can be conveniently captured by  enforcing a certain support pattern of the TT cores, yielding so-called  \emph{block-sparse} TTs.
This structure allows us to allocate resources much more efficiently---we can enlarge the dictionary while keeping the problem-relevant expressivity of the ansatz class and the computational resources required constant.

More formally, we define a \emph{degree map}  $w:[p] \rightarrow \mathbb{N}_0$, which assigns the \emph{degree} to the elements of a given function dictionary.
Without loss of generality we will assume that $w$ is a non-decreasing function.
The degree map provides us with a natural definition of the Laplace-like multivariate degree operator $L: \mathbb{R}^{p^d} \rightarrow \mathbb{R}^{p^d}$, which acts on tensor representations of multivariate functions,
\begin{equation}\label{eq:degree operator}
    L = \sum_{j=1}^d \id_p^{\otimes j-1} \otimes \Omega \otimes \id_p^{\otimes d-j}\,,
\end{equation}
where $\Omega = \diag(w(1),\dots,w(p))$ and $\id_p$ is the $p \times p$ identity matrix.
This operator is analogous to the bosonic particle number operator in quantum physics.

Suppose a tensor $\phi$ with a TT representation with cores $\{C_\ell\}_{\ell=1}^d$.
We define the \emph{left and right interface operators}
\begin{align}\label{eq:interface operators}
    L^{<\ell} &\coloneqq \sum_{j=1}^{\ell-1} \id^{\otimes j-1} \otimes \Omega \otimes \id^{\otimes \ell - j - 1},\\
    L^{>\ell} &\coloneqq \sum_{j=1}^{d-\ell} \id^{\otimes j-1} \Omega \otimes \id^{\otimes d-\ell-j},
\end{align}
and the \emph{left and right interface tensors}
\begin{equation}\label{eq:interface tensors}
    \begin{tikzpicture}[
        dot/.style = {circle,fill,minimum size=#1,inner sep=0pt,outer sep=0pt},
        dot/.default = 5pt
        ]
        \node[dot,label=above:$C_1$] at (0,0) (c1) {};
        \node[dot,label=above:$C_2$] [right=.5cm of c1] (c2) {};
        \node [right=.5cm of c2] (dots) {$\dots$};
        \node[dot,label=above:$C_{\ell-1}$] [right=.5cm of dots] (cd) {};
        \node [left=.5 of c1] (LHS) {$\phi^{<\ell} \coloneqq $};

        \node [below=.3cm of c1] (d1) {};
        \node [below=.3cm of c2] (d2) {};
        \node [below=.3cm of cd] (dd) {};
        \node [right=.3cm of cd] (right) {};

        \draw (d1) -- (c1) -- (dots.west);
        \draw (d2) -- (c2) -- (dots.west);
        \draw (dd) -- (cd) -- (dots.east);
        \draw (cd) -- (right);

        \node[dot,label=above:$C_{\ell+1}$] at (0,-1.2) (c1d) {};
        \node[dot,label=above:$C_{\ell+2}$] [right=.7cm of c1d] (c2d) {};
        \node [right=.5cm of c2d] (dotsd) {$\dots$};
        \node[dot,label=above:$C_{d}$] [right=.5cm of dotsd] (cdd) {};
        \node [left=.5 of c1d] (LHSd) {$\phi^{>\ell} \coloneqq $};

        \node [below=.3cm of c1d] (d1d) {};
        \node [below=.3cm of c2d] (d2d) {};
        \node [below=.3cm of cdd] (ddd) {};
        \node [left=.3cm of c1d] (leftd) {};

        \draw (d1d) -- (c1d) -- (dotsd.west);
        \draw (d2d) -- (c2d) -- (dotsd.west);
        \draw (ddd) -- (cdd) -- (dotsd.east);
        \draw (c1d) -- (leftd);
    \end{tikzpicture}
\end{equation}
for each interface at the $\ell$-th tensor core. 
When using matrix notation for the interface tensors, we will think of them as linear operators $\phi^{<\ell}: \mathbb{R}^{r_{\ell-1}}\rightarrow \mathbb{R}^{p^{\ell-1}}$ and $\phi^{>\ell}: \mathbb{R}^{p^{d-\ell}} \rightarrow \mathbb{R}^{r_\ell}$.

In ref.~\cite{gotte_block-sparse_2021-1}, the following theorem has been shown.

\begin{theorem}[Block-sparsity\label{thrm:block sparsity}]
    Suppose a TT $\phi$ with tensor cores $\{C_\ell\}_{\ell\in[d]}$ in \emph{left canonical form} with \emph{minimal ranks} $\{r_\ell\}_{\ell=1}^{d-1}$, such that 
    \begin{equation}
        L\phi = \lambda \phi,
    \end{equation}
    where $L$ is as in eq.~\eqref{eq:degree operator} and $\lambda \in \mathbb{N}_0$. 
    Then there exist a unitary gauge transformation $\{A_\ell\}_{\ell\in[d-1]}$ acting via eq.~\eqref{eq:gaugetrafo}, such that for each interface $\ell \in [d-1]$ the transformed interface tensors 
    satisfy
    \begin{align}
        \phi^{>\ell}L^{>\ell} &= \Lambda^{>\ell}\phi^{>\ell}\label{eq:block sparsity 1},\\
        L^{<\ell+1}\phi^{<\ell+1} &= \phi^{<\ell+1}\left(\lambda \id - \Lambda^{>\ell}\right)\label{eq:block sparsity 2}
    \end{align}
    for a set of diagonal matrices $\left\{\Lambda^{>\ell}\right\}_{\ell\in[d-1]}$ with non-increasing diagonal entries.
\end{theorem}

We provide a slightly simplified proof in \cref{app:block sparsity proof}.

Let us discuss the consequences of \cref{thrm:block sparsity} and see how it implies block-sparse structure of the tensor cores.
First, notice that eq.~\eqref{eq:block sparsity 1} is an eigenvalue equation, i.e. it states that the rows of $\phi^{>\ell}$ are left eigenvectors of $L^{>\ell}$ with eigenvalues given by the corresponding diagonal element of $\Lambda^{>\ell}$.
Similarly, eq.~\eqref{eq:block sparsity 2} implies that the rows of $\phi^{<\ell}$ are right eigenvectors of $L^{<\ell}$ with eigenvalues given by the corresponding diagonal elements of $\Lambda^{<\ell} \coloneqq \lambda \id - \Lambda^{>\ell-1}$.

Choosing any interface $\ell \in \{2,\dots,d-1\}$, we can decompose the eigenvalue equation as 
\begin{equation}
    \begin{tikzpicture}[
        dot/.style = {circle,fill,minimum size=#1,inner sep=0pt,outer sep=0pt},
        dot/.default = 5pt
        ]
        \node (lhs) at (0,0) {$\lambda \phi = L \phi =$};
        \node [dot,label=above:$\phi^{<\ell}$] (phi) [right=.5 of lhs] {};
        \node [dot,label=above:$C_\ell$] (cl) [right=.5 of phi] {};
        \node [dot,label=left:$L^{<\ell}$] (Ll) [below=.3 of phi] {};
        \node (Llb) [below=.3 of Ll]{};
        \node (id1) [below=.8 of cl] {};
        \node [dot,label=above:$\phi^{>\ell}$] (phiR) [right=.5 of cl] {};
        \node (idr1) [below=.8 of phiR] {};

        \draw[double] (Llb) -- (Ll) -- (phi);
        \draw (phi) -- (cl) -- (id1);
        \draw (cl) -- (phiR);
        \draw[double] (phiR) -- (idr1);
        
        \node (plus3) [right=.5 of phiR] {$+$};
        \node (plus) [below=1.1 of lhs] {$+$};

        \node [dot,label=above:$\phi^{<\ell}$] (phi2) [right=.5 of plus] {};
        \node [dot,label=above:$C_\ell$] (cl2) [right=.5 of phi2] {};
        \node [dot,label=left:$\Omega$] (omega) [below=.3 of cl2] {};
        \node (Llb2) [below=.3 of omega]{};
        \node (idl2) [below=.8 of phi2] {};
        \node [dot,label=above:$\phi^{>\ell}$] (phiR2) [right=.5 of cl2] {};
        \node (idr2) [below=.8 of phiR2] {};

        \draw[double] (idl2) -- (phi2);
        \draw (phi2) -- (cl2) -- (omega) -- (Llb2);
        \draw (cl2) -- (phiR2);
        \draw[double] (phiR2) -- (idr2);

        \node (plus2) [right=.5 of phiR2] {$+$};

        \node [dot,label=above:$\phi^{<\ell}$] (phi2) [right=.5 of plus2] {};
        \node [dot,label=above:$C_\ell$] (cl2) [right=.5 of phi2] {};
        \node (idl2) [below=.8 of phi2] {};
        \node [dot,label=above:$\phi^{>\ell}$] (phiR2) [right=.5 of cl2] {};
        \node [dot,label=right:$L^{>\ell}$] (lr) [below=.3 of phiR2] {};
        \node (Llb2) [below=.3 of lr]{};
        \node (idr2) [below=.8 of cl2] {};

        \draw[double] (idl2) -- (phi2);
        \draw (phi2) -- (cl2) -- (idr2);
        \draw (cl2) -- (phiR2);
        \draw[double] (phiR2) -- (lr) -- (Llb2);
        
        \node (eq) [below=3 of lhs] {$=$};
        \node [dot,label=above:$\phi^{<\ell}$] (phiL) [right=.3 of eq] {};
        \node [dot,label=above:$\Lambda^{<\ell}$] (L) [right=.6 of phiL] {};
        \node [dot,label=above:$C_\ell$] (c) [right=.6 of L] {};
        \node [dot,label=above:$\phi^{>\ell}$] (phiR) [right=.5 of c] {};
        \node (phiLB) [below=.3 of phiL] {};
        \node (phiRB) [below=.3 of phiR] {};
        \node (cB) [below=.3 of c] {};

        \draw[double] (phiLB) -- (phiL);
        \draw[double] (phiRB) -- (phiR);
        \draw (cB) -- (c);
        \draw (phiL) -- (L) -- (c) -- (phiR);

        \node (plus) [right=.3 of phiR] {$+$};
        \node [dot,label=above:$\phi^{<\ell}$] (phiL) [right=.3 of plus] {};
        \node [dot,label=above:$C_\ell$] (c) [right=.5 of phiL] {};
        \node [dot,label=above:$\Lambda^{>\ell}$] (L2) [right=.6 of c] {};
        \node [dot,label=above:$\phi^{>\ell}$] (phiR) [right=.6 of L2] {};
        \node (phiLB) [below=.3 of phiL] {};
        \node (phiRB) [below=.3 of phiR] {};
        \node (cB) [below=.3 of c] {};

        \draw[double] (phiLB) -- (phiL);
        \draw[double] (phiRB) -- (phiR);
        \draw (cB) -- (c);
        \draw (phiL) -- (c) -- (L2) -- (phiR);

        \node (plus) [below=.8 of L] {$+$};
        \node [dot,label=above:$\phi^{<\ell}$] (phiL) [right=.3 of plus] {};
        \node [dot,label=above:$C_\ell$] (c) [right=.5 of phiL] {};
        \node [dot,label=left:$\Omega$] (omega) [below=.3 of c] {};
        \node [dot,label=above:$\phi^{>\ell}$] (phiR) [right=.5 of c] {};
        \node (phiLB) [below=.8 of phiL] {};
        \node (phiRB) [below=.8 of phiR] {};
        \node (cB) [below=.3 of omega] {};
        \node (dot) [right=.5 of phiR] {,};

        \draw[double] (phiLB) -- (phiL);
        \draw[double] (phiRB) -- (phiR);
        \draw (cB) -- (c);
        \draw (phiL) -- (c) -- (phiR);
    \end{tikzpicture}
\end{equation}
with $\Omega = \diag(w(1),\dots,w(p))$, where the double line collects multiple indices into a single edge.
Since we assume that the TT ranks are minimal, $\phi^{<\ell}$ ($\phi^{>\ell}$) have full column (row) rank and the eigenvalue equation implies
\begin{equation}
    \begin{tikzpicture}[
        dot/.style = {circle,fill,minimum size=#1,inner sep=0pt,outer sep=0pt},
        dot/.default = 5pt
        ]
        \node (lmbd) {$\lambda$};
        \node[dot,label=above:$C_\ell$] (c) [right=.3 of lmbd] {};
        \node (left) [left=.3 of c] {};
        \node (right) [right=.3 of c] {};
        \node (below) [below=.3 of c] {};

        \draw (left) -- (c) -- (right);
        \draw (c) -- (below);

        \node (eq) [right=.3 of c] {$=$};
        \node[dot,label=above:$C_\ell$] (c) [right=1.2 of eq]{};
        \node (right) [right=.3 of c] {};
        \node[dot,label=above:$\Lambda^{<\ell}$] (lambda) [left=.5 of c] {};
        \node (left) [left=.3 of lambda] {};
        \node (below) [below=.3 of c] {};

        \draw (left) -- (c) -- (right);
        \draw (c) -- (below);

        \node (plus) [right=-.2 of right] {$+$};
        \node[dot,label=above:$C_\ell$] (c) [right=.4 of plus] {};
        \node[dot,label=above:$\Lambda^{>\ell}$] (lambda) [right=.5 of c] {};
        \node (right) [right=.3 of lambda] {};
        \node (left) [left=.3 of c] {};
        \node (below) [below=.3 of c] {};

        \draw (left) -- (c) -- (lambda) -- (right);
        \draw (c) -- (below);

        \node (eq) [right=-.2 of right] {$+$};
        \node[dot,label=above:$C_\ell$] (c) [right=.3 of eq] {};
        \node (left) [left=.3 of c] {};
        \node (right) [right=.3 of c] {};
        \node[dot,label=left:$\Omega$] (omega) [below=.2 of c] {};
        \node (below) [below=.2 of omega] {};
        \node (dot) [right=-.2 of right] {.};

        \draw (left) -- (c) -- (right);
        \draw (c) -- (omega) -- (below);
    \end{tikzpicture}
\end{equation}
If we fix the remaining physical index to $i\in[p]$, we obtain the matrix equation
\begin{equation}\label{eq:interpreting block sparsity}
    (\lambda - w(i)) \left(C_\ell\right)_i = \Lambda^{<\ell} \left(C_\ell\right)_i + \left(C_\ell\right)_i \Lambda^{>\ell}.
\end{equation}
Therefore, for all $i\in[p]$, we can assign each block-row of $\left(C_\ell\right)_i$ to an eigenvalue of $L^{<\ell}$ by which it is multiplied in the first term on the RHS of eq.~\eqref{eq:interpreting block sparsity}.
Similarly we can assign each block-column of $\left(C_\ell\right)_i$ to an eigenvalue of $L^{>\ell}$ by which it is multiplied in the second term on the RHS of eq.~\eqref{eq:interpreting block sparsity}.
Now, eq.~\eqref{eq:interpreting block sparsity} tells us that blocks of $\left(C_\ell\right)_i$, which correspond to eigenvalues $\lambda^{<\ell}$, $\lambda^{>\ell}$ of $L^{<\ell}$ and $L^{>\ell}$ respectively, can be non-zero only if 
\begin{equation}\label{eq:block non-zero cond}
    \lambda^{>\ell} + \lambda^{<\ell} + w(i) = \lambda\, .
\end{equation}
Hence, each non-zero block of $C_\ell$ connects eigenvectors of $L^{<\ell}$ with eigenvectors of $L^{>\ell}$ to fulfil the eigenvalue equation $L\phi = \lambda\phi$.
In this way \cref{thrm:block sparsity} implies block-sparsity of the tensor cores of an eigenvector of $L$.

In order to analyze the maximum block sizes, consider a block of $C_\ell$ that corresponds to the eigenvalues $\lambda^{<\ell}$ and $\lambda^{>\ell}$.
The block size is limited by the number of times each eigenvalue appears in the spectrum of $L^{<\ell}$ and $L^{>\ell}$, respectively.
The interface operators are diagonal, with elements given by $\sum_{\ell'=1}^{\ell-1} w(i_{\ell'})$ and $\sum_{\ell'=\ell+1}^d w(j_{\ell'})$, respectively, so the maximum size $s_\ell \times t_\ell$ of the considered block, is the number of solutions $(i_1,\dots,i_{\ell-1}) \in [p]^{\ell-1}$ and $(j_{\ell+1},\dots,j_d) \in [p]^{d-\ell}$ to
\begin{equation}\label{eq:block sizes}
    \sum_{\ell'=1}^{\ell-1} w(i_{\ell'}) = \lambda^{<\ell} \quad \text{and} \quad \sum_{\ell'=\ell+1}^{d} w(i_{\ell'}) = \lambda^{>\ell},
\end{equation}
respectively.
Note that if we allowed $w$ to take negative values, the maximum block sizes and number of blocks would become very large.

To obtain a \emph{low-rank} block-sparse TT, we enforce the block-sparse structure and limit the block sizes to some maximum value $\rho$.

\paragraph*{Example 1: Monomial dictionary.}
Consider the monomial dictionary
\begin{equation}\label{eq:monomial dict}
    \Psi_i(x) = x^{i-1} \ , \quad i\in[p]
\end{equation}
and the polynomial degree function $w: i \mapsto i-1$, so that $\image(w) = \{0,\dots,p-1\}$.
For physical index $i$ at the $\ell$-th mode, there are $\lambda - i + 2$ solutions $\lambda^{<\ell}$ and $\lambda^{>\ell}$ to eq.~\eqref{eq:block non-zero cond}, so the matrix $\left(C_\ell\right)_i$ has $\lambda - i + 2$ non-zero blocks.
A combinatorial argument shows that the number of solutions to eq.~\eqref{eq:block sizes} is
\begin{equation}
    s_\ell = \genfrac(){0pt}{0}{\lambda^{<\ell} + \ell - 1}{\ell - 2} \ , \quad t_\ell = \genfrac(){0pt}{0}{\lambda^{>\ell} + d - \ell}{\ell - 2},
\end{equation}
which gives us the maximum block sizes.
For concreteness, choose $\lambda = 3$ and $p=3$.
Now the block-sparse structure becomes
\begin{equation}\label{eq:monomial block sparsity}
    \begin{split}
        \left(C_\ell\right)_1 = \begin{pmatrix} *&0&0&0\\0&*&0&0\\0&0&*&0\\0&0&0&* \end{pmatrix}&,\quad
        \left(C_\ell\right)_2 = \begin{pmatrix} 0&*&0&0\\0&0&*&0\\0&0&0&*\\0&0&0&0 \end{pmatrix},\\
        &\hspace{-1.6cm}\left(C_\ell\right)_3 = \begin{pmatrix} 0&0&*&0\\0&0&0&*\\0&0&0&0\\0&0&0&0 \end{pmatrix},
    \end{split}
\end{equation}
where $*$ indicates the non-zero blocks. 
Here the $i$-th row corresponds to $\lambda^{<\ell}=i-1$ and the $j$-th column to $\lambda^{>\ell} = 4 - j$.

\paragraph*{Example 2: Trigonometric dictionary.}
Consider the trigonometric dictionary 
\begin{equation}\label{eq:trig dict}
    \{1,\sin x,\cos x\}
\end{equation}
and $w(1)=0,w(2)=1,w(3)=1$, which counts the number of sines and cosines.
Choose $\lambda = 3$ and $\ell \in \{4,\dots,d-4\}$.
The block-sparse structure now becomes 
\begin{equation}\label{eq:trig block sparsity}
    \begin{split}
        \left(C_\ell\right)_1 = \begin{pmatrix} *&0&0&0\\0&*&0&0\\0&0&*&0\\0&0&0&* \end{pmatrix}& \ , \quad
        \left(C_\ell\right)_2 = \begin{pmatrix} 0&*&0&0\\0&0&*&0\\0&0&0&*\\0&0&0&0 \end{pmatrix},\\
        \left(C_\ell\right)_3 &= \begin{pmatrix} 0&*&0&0\\0&0&*&0\\0&0&0&*\\0&0&0&0 \end{pmatrix},
    \end{split}
\end{equation}
where the $i$-th row corresponds to $\lambda^{<\ell} =i-1$ and the $j$-th column to $\lambda^{>\ell} = 4 - j$.

\paragraph*{Limiting the degree.}
\cref{thrm:block sparsity} states that \emph{fixed degree functions} admit a \emph{block-sparse TT representation}.
An often more natural ansatz class are functions with bounded degree.
Fortunately, such functions also admit a block-sparse TT, as can be seen from the following argument.

Suppose a function $f:\mathbb{R}^d$ with a bounded degree $\lambda$.
Hence, it can be written as a linear combination 
\begin{equation}
    f = \sum_j f_j,
\end{equation}
where each $f_j$ has a fixed degree $\le \lambda$ and, hence, admits a block-sparse decomposition.
The number of terms $n(\lambda)$ in this sum is bounded by the number of options for a degree $\le \lambda$.
Suppose that each $f_j$ has a block-sparse TT decomposition $\{C_\ell^{(j)}\}_{\ell\in[d]}$.
Now we can represent the function $f$ by a TT with tensor cores $\{C_\ell\}_{\ell \in [d]}$, where for each $\ell \in\{2,\dots,d-1\}$ and a value $i\in[p]$ of the physical index, the corresponding matrix $\left(C_\ell\right)_i$ is block-diagonal with blocks $\left(C_\ell^{(j)}\right)_i$, the tensor core $C_1$ is given by concatenating $C_1^{(j)}$ and the tensor core $C_d$ is
\begin{equation}
    \begin{tikzpicture}[
        dot/.style={circle,fill,minimum size=#1,inner sep=0pt,outer sep=0pt},
        dot/.default=5pt
        ]
        \node [dot,label=above:$C_d$] (cd) {};
        \node (left) [left=.3 of cd] {};
        \node (right) [right=.3 of cd] {$j$};
        \node (bottom) [below=.3 of cd] {};
        \node (eq) [right=.3 of right] {$=$};
        \node [dot,label=above:$\tilde C_d^{(j)}$] (cdR) [right=.8 of eq] {};
        \node (bottom2) [below=.5 of cdR] {};
        \node (left2) [left=.5 of cdR] {};
        \node (dot) [right=.2 of cdR] {,};

        \draw (left) -- (right);
        \draw (bottom) -- (cd);
        \draw (left2) -- (cdR) -- (bottom2);
    \end{tikzpicture}
\end{equation}
where $C_d^{(j)}$ has $C_d^{(j)}$ as its $j$-th block, with all other blocks zero.
We sum over the blocks corresponding to each $f_j$ by contracting the right index of $C_d$ with a vector of all ones $\mathbf{1}_{n(\lambda)}$.
Hence, the TT representation of $f$ becomes
\begin{equation}\label{eq:degree limited TT}
    \begin{tikzpicture}[
        dot/.style={circle,fill,minimum size=#1,inner sep=0pt,outer sep=0pt},
        dot/.default=5pt
        ]
        \node [dot,label=above:$C_1$] (c1) {};
        \node (c1R) [right=.3 of c1] {};
        \node (c1B) [below=.5 of c1] {};
        \node (dots) [right=.5 of c1] {$\dots$};
        \node [dot,label=above:$C_d$] (cd) [right=.5 of dots] {};
        \node (cdL) [left=.3 of cd] {};
        \node (cdB) [below=.5 of cd] {};
        \node [dot,label=right:$\id_{n(\lambda)}$] (id) [right=.5 of cd] {};
        \node [dot,label=right:$\mathbf{1}_{n(\lambda)}$] (one) [below=.5 of id] {};
        \node (dot) [right=1 of id] {,};
        
        \draw (c1B) -- (c1) -- (c1R);
        \draw (cdL) -- (cd) -- (id) -- (one);
        \draw (cd) -- (cdB);
    \end{tikzpicture}
\end{equation}
which is a block-sparse TT.

\section{Efficient TT representation of dynamical laws}\label{sec:efficient}

In this section we will show how generic properties of dynamical systems imply efficient TT representations of their dynamical laws.
We write the function $f: \mathcal{S} \rightarrow \mathbb{R}^d$ in eq.~\eqref{eq:dynamical law} using the decomposition \eqref{eq:product basis decomp} as 
\begin{equation}
    f_k(x) = \sum_{i_1,\dots,i_d} \Theta_{k;i_1,\dots,i_d} \Phi_{i_1,\dots,i_d}(x).
\end{equation}
We call the TT decomposition of $\Theta_k$ the TT representation of $f_k$.
Combining the tensors for $k\in[d]$, we obtain the tensor $\Theta \in \mathbb{R}^{d \times p^d}$, which can be written as a single TT via
\begin{equation}\label{eq:total TT}
\begin{tikzpicture}
    \node at (-1.4,0) {$\Theta_{k;i_1,\ldots ,i_d} =$};
    \node [above] at (0,0) {$C_{1}$};
    \draw (0,-0.5) -- (0,0) -- (1,0);
    \draw [fill] (0,0) circle (0.5ex);
    \node [below] at (0,-0.5) {$i_1$};
    
    \node [above] at (1,0) {$C_{2}$};
    \draw (1,-0.5) -- (1,0) -- (1.5,0);
    \draw [fill] (1,0) circle (0.5ex);
    \node [below] at (1,-0.5) {$i_2$};

    \node at (2,0) {$\dots$};
    
    \node [above] at (3,0) {$C_{d-1}$};
    \draw (2.5,0) -- (4,0) -- (4,-0.5);
    \draw (3,0) -- (3,-0.5);
    \draw [fill] (3,0) circle (0.5ex);
    \node [below] at (3,-0.5) {$i_{d-1}$};

    \node [above] at (4,0) {$C_{d}$};
    \draw [fill] (4,0) circle (0.5ex);
    \node [below] at (4,-0.5) {$i_d$};
    \draw (4,0) -- (4.5,0);
    \node [right] at (4.5,0) {$k$};

    \node [below] at (5,0) {,};
\end{tikzpicture}
\end{equation}
providing a TT representation of the dynamical law.

In order to show that a function admits an efficient representation, we will bound its \emph{separation rank}, defined as follows.

\begin{definition}[Separation rank]
    Suppose a multivariate function $h(x):\mathbb{R}^d \rightarrow \mathbb{R}$.
    We say that $h(x)$ has \emph{separation rank $s$ with respect to a bipartition $\mathcal{P}_k \coloneqq (\{x_1,\dots,x_k\},\{x_{k+1},\dots,x_d\}) \eqqcolon (\mathcal{P}^\text{left}_k,\mathcal{P}^\text{right}_k)$}, if the smallest set of functions $\{g_\ell^\text{left}(\mathcal{P}_k^\text{left}), g_\ell^\text{right}(\mathcal{P}^\text{right}_k)\}_{\ell\in[\tilde s]}$, such that 
    \begin{equation}
        h(x) = \sum_{\ell \in [\tilde s]} g_\ell^\text{left}(\mathcal{P}_k^\text{left}) \otimes g_\ell^\text{right}(\mathcal{P}_k^\text{right}),
    \end{equation}
    has $\tilde s = s$.
\end{definition}

It can be shown that the minimal ranks $r_k$ of a TT representation of a function $h(x)$ are equal to the separation ranks with respect to $\mathcal{P}_k$.
For a formal proof of this statement, see ref.~\cite{holtz_manifolds_2012}.

The first generic property of dynamical systems that has been used in \cite{goesmann_tensor_2020} to show efficient TT representations of dynamical laws is \emph{locality} in one-dimension.

\begin{definition}[One dimensional interacting system, Definition 4 in ref.~\cite{goesmann_tensor_2020}]\label{def:one-dim interacting}
    A dynamical system, governed by the dynamical law $f:\mathcal{S} \subset \mathbb{R}^d$, is \emph{one-dimensional with interaction length $L$ and separation rank $N$}, if there exists a function set $\{g_i(x)\}_{i=1}^{\tilde p}$ and for each $k \in [d]$ an index set $\mathcal{I}_k \subset [\tilde p]^{2L + 1}$ with $|\mathcal{I}_k|\le N$, such that 
    \begin{equation}
        f_k(x) = \sum_{i_{k-L},\dots,i_{k+L} \in \mathcal{I}_k} g_{i_{k-L}}(x_{k - L}) \dots g_{i_{k+L}}(x_{k+L}),
    \end{equation}
    where we set $g_i = 1$ for $i \le 0$ and $i \ge d+1$.
\end{definition}

By bounding the separation ranks of the dynamical laws of one dimensional interacting systems with respect to bipartitions $\mathcal{P}_k$, the authors prove the following theorem.

\begin{theorem}[Efficient TT decomposition of one-dimensional interacting systems, Theorem 5 in ref.~\cite{goesmann_tensor_2020}]\label{thrm:TT rep of local systems}
    Suppose a one dimensional interacting system with interaction length $L$ and separation rank $N$, governed by the dynamical law $f(x): \mathcal{S} \subset \mathbb{R}^{p^d} \rightarrow \mathbb{R}^d$.
    Then we have  the following.
    \begin{enumerate}
        \item Each $f_k(x)$ admits a TT representation with rank $r \le N$ and
        \item the function $f$ admits a TT representation with ranks $r_k \le k - L + 1 + 2NL$.
    \end{enumerate}
\end{theorem}

For many systems of interest exact locality as in \cref{def:one-dim interacting} is too strong of an assumption.
We would like to be able to use TTs also for systems that don't have a sharp bound on the interaction length, but where instead the interactions decay with distance. 
This is formalized by the following definition.

\begin{definition}[One-dimensional systems with algebraically decaying interactions]\label{def:decaying interactions}
    A dynamical system, governed by the dynamical law $f:[0,1)^d \rightarrow \mathbb{R}^d$, is \emph{one-dimensional with $(\chi,g)$-algebraically decaying interactions and separation rank $N$}, if there exists a function set $\{g_{k,L}\}_{k,L\in[d]}$, where for all $k,L \in [d]$ the function $g_{k,L}$ depends non-trivially only on $x_i$ for $i \in [k-L,k+L]$, satisfies $\|g_{k,L}\|_2 \le g$ and has separation rank bounded by $N$ with respect to bipartitions $\mathcal{P}_\ell = \{\{x_{k-L},\dots,x_\ell\},\{x_{\ell+1},\dots,x_{k+L}\}\}$ for all $\ell \in \{k-L,\dots,d+L\}$, such that 
    \begin{equation}\label{eq:algebraically decaying}
        f_k(x) = \sum_{L \in [d]} L^{-\chi} g_{k,L}(x).
    \end{equation}
\end{definition}

Note that in this definition we demand that the state space $\mathcal{S} = [0,1)^d$. 
This is important so that $\|g_L\|_2$ does not change with $d$.
For these systems, we can show the following theorem.

\begin{theorem}[Approximate locality]\label{thrm:approximate locality}
    Suppose a one-dimensional system with $(\chi,g)$-algebraically decaying interactions and separation rank $N$, governed by the dynamical law $f:[0,1)^d \rightarrow \mathbb{R}^d$, that can be written as
    \begin{equation}
        f_k(x) = \sum_{L \in [d]} L^{-\chi} g_{k,L}(x),
    \end{equation}
    where $\{g_{k,L}\}_{k,L\in[d]}$ satisfies the assumptions of \cref{def:decaying interactions}.
    If $\chi >1$, then for any $\tilde L \in [d]$ there exists a one-dimensional interacting system with interaction length $\tilde L$ and separation rank $N\tilde L$, such that 
    \begin{equation}
        \|f_k(x) - \tilde f_k^{\tilde L}(x)\|_2 \le c_1(\chi,\tilde L) g,
    \end{equation}
    where 
    \begin{equation}
        c_1(\chi,\tilde L) = (\tilde L + 1)^{-\chi} \frac{\tilde L + \chi}{\chi - 1} = \mathcal{O}(\tilde L^{1-\chi}).
    \end{equation}
\end{theorem}

The proof is given in \cref{app:approximate locality proof}.

\cref{thrm:approximate locality} allows us to approximate systems with algebraically decaying interactions by strictly local systems with bounded error, which is independent of $d$.\footnote{Note that if we were instead looking at the $L_2$ error of the full function $\|f - \tilde f\|_2$, we would get an additional factor of $d$.}
This is formalized in the following corollary.

\begin{corollary}[Low rank TTs for algebraically decaying interactions]\label{cor:approx locality}
    Suppose a system with $(\chi,g)$-algebraically decaying interactions with separation rank $N$ governed by the dynamical law $f(x):\mathcal{S}\subset \mathbb{R}^d \rightarrow \mathbb{R}^d$.
    Furthermore, suppose that $\chi>1$.
    Then there exists an $\varepsilon$-approximate TT representation of each $f_k(x)$ with rank
    \begin{equation}
        r \le N\left[\left(\frac\chi{\chi-1}\frac g\varepsilon \right)^{\frac1{\chi-1}} - 1\right].
    \end{equation}
\end{corollary}

The details of the proof are given in \cref{app:approx locality corollary proof}.

Many systems of physical interest are \emph{$K$-mode interacting systems}, formalized by the following definition.

\begin{definition}[$K$-mode interacting systems]\label{def:k-mode interacting}
    A dynamical system, governed by the dynamical law $f:\mathcal{S} \subset \mathbb{R}^d$, is \emph{$K$-mode interacting with separation rank $N$}, if there exists a function set $\{g_i(x)\}_{i=1}^{\tilde p}$, for each $k \in [d]$ a constant $\mathcal{K}_k$, such that for each $\ell \in [\mathcal{K}_k]$ there are distinct subsets $\mathcal{J}_\ell^{(k)} = \{j_1,\dots,j_K\} \subset [d]$ with $\mathcal{J}_{\ell}^{(k)} \ni k$ and subsets $\mathcal{I}_{\ell}^{(k)}\subset[\tilde p]^K$ with $|\mathcal{I}_\ell^{(k)}| \le N$, such that 
    \begin{equation}
        f_k(x) = \sum_{\ell \in [\mathcal{K}_k]} \sum_{(i_{j_1},\dots,i_{j_K}) \in \mathcal{I}_{\ell}^{(k)}} g_{i_{j_1}}(x_{j_1}) \dots g_{i_{j_K}}(x_{j_K}).
    \end{equation}
\end{definition}

An example of a $2$-mode interacting system is, e.g., a collection of gravitationally interacting particles, where the total force on each particle is given by the sum of pair-wise forces with respect to the remaining particles. 
We will now show that if a system is $K$-mode interacting with $K \lessapprox 5$, it admits an efficient TT representation.

\begin{theorem}[Efficient TT representation of $K$-mode interacting systems]\label{thrm:TT rep of K-mode sys}
    Suppose a $K$-mode interacting system with separation rank $N$, governed by a dynamical law $f(x):\mathcal{S} \subset \mathbb{R}^{p^d} \rightarrow \mathbb{R}^d$. Then 
    \begin{enumerate}
        \item each $f_k$ admits a TT representation with rank $r \le N \genfrac(){0pt}{0}{d-1}{K-1} = N\mathcal{O}(d^{K-1})$ and
        \item the function $f$ admits a TT representation with ranks 
        \begin{equation}
            r_k \le c_2(N,d,k) + k \genfrac(){0pt}{0}{k-1}{K-1} + 1,
        \end{equation}
        where $c_2(N,d,k) = \mathcal{O}(Nd^K)$ is a combinatorial factor defined in \cref{app:K-mode proof}.
    \end{enumerate}
\end{theorem}

To prove \cref{thrm:TT rep of K-mode sys}, we adapt the techniques of the proof of \cref{thrm:TT rep of local systems} from ref.~\cite{goesmann_tensor_2020}, see \cref{app:K-mode proof} for details.

Already for $K \gtrapprox 5$, although $c_2(N,d,k)$ is polynomial in $d$, the polynomial degree makes working with such TTs prohibitively expensive even for modest $d$. 
However, many systems of physical interest are known to have $K=2,3,4$ and hence admit an efficient TT representation.
It is important to note that this result does not rely on the underlying systems being one-dimensional.
Finally, note that if a dynamical system has \emph{at most} $K$-mode interactions, it still admits an efficient TT decomposition, since the rank of a TT is sub-additive.

General conditions for the approximability of multivariate functions by a low rank TT in terms of tail control of the singular value spectrum of the matrix unfoldings of their coefficient tensors are derived in Ref.~\cite{bachmayr_tensor_2016}.
As we discuss in \cref{app:conditions_on_lowrankness}, these results are in a precise sense analog to the control of matrix product state approximations of quantum states based on entropy scaling conditions derived in the quantum many-body literature \cite{schuch_entropy_2008,verstraete_matrix_2006}.

\section{Gauge mediated weight sharing (ALS-GMWS)}\label{sec:gauge mediated weight sharing}

\subsection{Self-similarity}\label{ssec:self-similarity}

Additional structure in the system, known prior to learning, can cause certain modes to play the same role in dynamical laws for multiple modes, implying that we would like the corresponding tensor cores to be equal.
We call this \emph{self-similarity}.
The different roles each mode can play are referred to as \emph{activation types}.
For example, in the case of a one-dimensional dynamical system with interaction length $L$, \cref{def:one-dim interacting}, the $j$-th mode plays the same role in all functions $f_k$ with $k < j-L$, namely that the mode is to the left and outside of the interaction range.
Similarly the role of the mode is the same for all $k > j+L$.
Hence, such systems have $2L + 3$ activation types.

Self-similar systems with $\alpha$ activation types can be described by a set of $d\alpha$ tensor cores $\{C_\ell^{(j)}\}_{\ell \in [d], j \in [\alpha]}$.
The recipe to build the corresponding tensor train representations of the dynamical laws $f_k(x)$ can be encoded in a selection table $S \in [\alpha]^{d\times d}$ as in 
\begin{equation}
    \begin{tikzpicture}[
        dot/.style={circle,fill,minimum size=#1,inner sep=0pt,outer sep=0pt},
        dot/.default=5pt
        ]
        \node [dot,label=above:$C_1^{(S_{k,1})}$] (c1) {};
        \node (c1R) [right=.3 of c1] {};
        \node (c1B) [below=.5 of c1] {};
        \node (dots) [right=.5 of c1] {$\dots$};
        \node [dot,label=above:$C_d^{(S_{k,d})}$] (cd) [right=.5 of dots] {};
        \node (cdL) [left=.3 of cd] {};
        \node (cdB) [below=.5 of cd] {};
        \node [dot,label=right:$\id_{n(\lambda)}$] (id) [right=.5 of cd] {};
        \node [dot,label=right:$\mathbf{1}_{n(\lambda)}$] (one) [below=.5 of id] {};
        \node (dot) [right=1 of id] {,};
        
        \draw (c1B) -- (c1) -- (c1R);
        \draw (cdL) -- (cd) -- (id) -- (one);
        \draw (cd) -- (cdB);
    \end{tikzpicture}
\end{equation}
where block-sparse representation of functions with bounded degree eq.~\eqref{eq:degree limited TT} is used.
We call such systems $S$-self-similar.

In the case of one-dimensional interacting systems with interaction length $L$, the selection table takes the form
\begin{equation}\label{eq:local selection table}
    S_{ij}^L = 
    \begin{cases}
        1 & j < i-L,\\
        2L + 3 & j > i+L,\\
        j-i+L+2 & \text{otherwise.}
    \end{cases}
\end{equation}

\subsection{ALS optimization of self-similar systems}\label{ssec:ALS}

In learning dynamical laws, given data $(x_i,y_i) \in \mathbb{R}^d \times \mathbb{R}^d$ for $i\in[M]$ with $y_i \approx f(x_i) \ \forall i$, we would like to identify an element $\hat f$ of some ansatz class that minimizes the empirical loss
\begin{equation}
    \mathcal{L}_\text{emp}(\hat f) \coloneqq \sum_{i\in[M]} \left\|y_i - \hat f(x_i)\right\|_2^2.
\end{equation}
In the previous sections we have shown that a natural choice of an ansatz class for many dynamical systems of interest are block-sparse low rank TTs with self-similarity given by the selection table eq.~\eqref{eq:local selection table}.

Previously, in ref.~\cite{goesmann_tensor_2020}, \emph{alternating least squares} (ALS)  optimization (and a rank-adaptive variant \cite{grasedyck_stable_2019}) has been used to minimize $\mathcal{L}_\text{emp}$ over low rank TTs.
In the ALS procedure, the tensor cores are iterated over in sweeps, at each step solving a linear least squares problem to minimize the empirical loss as a function of the given core, until convergence.
The ALS algorithm can be adapted to block-sparse tensor trains by restricting each contraction in the algorithm to indices labeling elements that are non-zero in the block-sparse structure \cite{gotte_block-sparse_2021-1}.

However, optimization over systems with self-similarity requires more care.
In ref.~\cite{goesmann_tensor_2020} a selection tensor approach has been taken, where the dynamical law is written as
\begin{equation}\label{eq:selection tensor}
    \begin{tikzpicture}[
        dot/.style={circle,fill,minimum size=#1,inner sep=0pt,outer sep=0pt},
        dot/.default=5pt
        ]
        \node [dot,label={[label distance=-5]135:$C_1$}] (c1) {};
        \node [dot,label={[label distance=-5]135:$C_2$}] [right=.6 of c1] (c2) {};
        \node [right=.3 of c2] (dots) {$\dots$};
        \node [dot,label={[label distance=-5]135:$C_d$}] [right=.3 of dots] (cd) {};
        \node [dot,label={[label distance=-5]135:$\mathcal{T}$}] [above=.3 of dots] (selection) {};
        \node [above=.3 of selection] (above) {$k$};
        \node [below=.3 of c1] (c1b) {};
        \node [below=.3 of cd] (cdb) {};
        \node [below=.3 of c2] (c2b) {};
        \node (lhs) [left=.5 of c1] {$\theta_k = $};
        \node (dot) [right=.1 of cd] {,};

        \draw (c1) to[in=180,out=90] (selection) to[in=90,out=0] (cd);
        \draw (c2) to[in=190,out=90] (selection);
        \draw (selection) -- (dots);
        \draw (c1) -- (c2) -- (dots) -- (cd);
        \draw (c1b) -- (c1);
        \draw (cdb) -- (cd);
        \draw (c2b) -- (c2);
        \draw (selection) -- (above);
    \end{tikzpicture}
\end{equation}
where $\theta_k$ is the TT representation of $f_k$ and $\mathcal{T}$ is the selection tensor that, given $k\in[d]$, picks the correct activation type.
In this form, block-sparse ALS can be directly applied.
However, this method leads to a numerical instabilities in the optimization, limiting its scalability to $d\lessapprox20$.

The core problem is that the unitary gauge freedom is not properly taken care of.
To see this, suppose that the ALS sweeps are performed from left to right and that the tensor core $C_\ell^{(j)}$ corresponding to the mode $\ell\in[d]$ and activation type $j\in[\alpha]$ is being optimized. 
The left neighbour of this tensor core can have a different activation type for different choices of $k$.
Since each of the left neighbour activation types has been optimized in a separate ALS step, they are each written in a different gauge.
Using the form eq.~\eqref{eq:selection tensor} does not 
allow the gauge differences to be corrected for and therefore we rely on the gauges to converge as more ALS sweeps are performed.

We will now introduce a new algorithm for optimization of block-sparse low rank TTs with self-similarity, called \emph{ALS with gauge mediated weight sharing} (ALS-GMWS), that allows the left neighbour gauge to be adjusted in each step.
It will be useful to define \emph{restricted empirical cost}, which, given a subset $E \subset [d]$, is given by
\begin{equation}
    \mathcal{L}^E_\text{emp}(\hat f) \coloneqq \sum_{i\in[M]} \left\|\left(y_i\right)_E - \left(\hat f(x_i)\right)_E \right\|_2^2,
\end{equation}
where the notation $\left(v\right)_E$ for $v \in \mathbb{R}^d$ denotes the restriction of $v$ to the subspace defined by the set of index values $E$. 

Assume without loss of generality that the ALS sweeps are performed left to right.
Suppose we want to optimize the core $C^{(j)}_\ell$, for $\ell \in [d]$ and $j \in [\alpha]$.
Let $\mathcal{E}\subset [d]$ be the set of indices of functions where this core is used, as determined by the selection table $S$.
If $\ell=1$, we don't need to worry about the gauge at all, since we are optimizing left to right.
We simply perform a block-sparse ALS step to find $C_\ell^{(j)}$ that minimizes $\mathcal{L}^\mathcal{E}_\text{emp}(\hat f)$ with all the other tensor cores fixed.

Now consider $\ell\ge2$.
Among the functions labeled by $e \in \mathcal{E}$, we are only sure that the gauge is the same in the equivalence classes defined by the values of $S_{e,\ell-1}$, or, in other words, only if the left neighbour, the tensor core of the $\ell-1$-th mode, has the same activation type.
Hence, we divide $\mathcal{E}$ into disjoint sets $\{E_a\subset \mathcal{E}\}_{a\in[\alpha]}$ labeled by the activation type of $\ell-1$-th mode, such that, for each $a\in[\alpha]$, $S_{e,\ell-1} = a$ for all $e \in E_a$.
In order to use the most information available, we find $\tilde a \in [\alpha]$, such that $|E_{\tilde a}| \ge |E_a| \ \forall a \in [\alpha]$, and perform the block-sparse ALS step to find the $C^{(j)}_\ell$ that minimizes $\mathcal{L}^{E_{\tilde a}}_\text{emp}(\hat f)$, with all the other tensor cores fixed.

Now that we found the optimal $C^{(j)}_\ell$ with respect to $E_{\tilde a}$, we optimize the gauge of $C_{\ell-1}^{(a)}$ for $a \neq \tilde a$, so that using the newly found core in the corresponding functions is justified.
Hence, for each $a \in [\alpha]$, such that $a\neq \tilde a$ and $|E_a| \neq 0$, we want to find a gauge fixing unitary $U_a$ that minimizes $\mathcal{L}^{E_a}_\text{emp}(\hat f)$ under the transformation
\begin{equation}\label{eq:gauge fixing}
    \begin{tikzpicture}[
        dot/.style={circle,fill,minimum size=5pt,inner sep=0pt,outer sep=0pt}
        ]
        \node [dot,label=above:$C^{(a)}_{\ell-1}$] (c) {};
        \node (left) [left=.3 of c] {};
        \node (right) [right=.3 of c] {};
        \node (below) [below=.3 of c] {};

        \node (arrow) [right=.7 of c] {$\mapsto$};

        \node [dot,label=above:$C^{(a)}_{\ell-1}$] (c2) [right=.5 of arrow] {};
        \node (left2) [left=.3 of c2] {};
        \node (below2) [below=.3 of c2] {};
        \node [dot,label=above:$U_a$] (u) [right=.5 of c2] {};
        \node (right2) [right=.3 of u] {};
        \node (dot) [right=.5 of u] {,};

        \draw (left) -- (c) -- (right);
        \draw (below) -- (c);
        \draw (left2) -- (c2) -- (u) -- (right2);
        \draw (below2) -- (c2);
    \end{tikzpicture}
\end{equation}
with all other tensor cores fixed.
In order to preserve the block-sparse structure, $U_a$ needs to be block-diagonal with block sizes given by the block-column widths of $C^{(a)}_{\ell-1}$ and block-row heights of $C^{(j)}_\ell$.
In practice, we optimize over all block-diagonal matrices, although optimization over block-diagonal unitaries is possible and could lead to improvements.

The ALS-GMWS algorithm is summarized in Algorithm \ref{alg:als-gmws}.

\begin{algorithm}[t!]
\SetKwInOut{Input}{input}\SetKwInOut{Output}{output}
\Input{Data pairs $(x_i,y_i)\in\mathbb{R}^d\times \mathbb{R}^d$, $i=1,\ldots,M$, number of activation types $\alpha$, selection table $S$, maximum block size $\rho$}
\Output{$d\alpha$ block-sparse tensor train components.}
\BlankLine
For $\ell=[d], \ j=[\alpha]$ initialize block-sparse $C^{(j)}_{\ell}$ with maximum block size $\rho$\;
\While{not converged}{
    \For{$\ell = 1,\dots,d$}{
        \For{$k = 1,\dots,\alpha$}{
            Let $\mathcal{E} = \{e\in[d]: S_{e,\ell}=j\}$\;
            \If{$\ell=1$ \texttt{or} $|\mathcal{E}| = 1$}{
                Find $C^{(k)}_\ell$ that minimizes $\mathcal{L}^\mathcal{E}_\text{emp}$ with all other tensor cores fixed\;
            }
            \Else{
                For each $a\in[\alpha]$ let $E_a = \{e\in\mathcal{E}: S_{e,\ell-1} = a\}$\;
                Choose $\tilde a \in [\alpha]$, s.t. $|E_{\tilde a}| \ge |E_a| \ \forall a\in[\alpha]$\;
                Find $C^{(k)}_\ell$ that minimizes $\mathcal{L}^{E_{\tilde a}}_\text{emp}$ with all other tensor cores fixed\;
                \For{$a \in [\alpha]: \ a\neq \tilde a, \ |E_a| \neq 0$}{
                    Find a block-diagonal unitary $U$ that minimizes $\mathcal{L}^{E_a}_\text{emp}$ after the transformation eq.~\eqref{eq:gauge fixing} with all the other tensor cores fixed\;
                    Set $C^{(a)}_{\ell-1} \leftarrow C^{(a)}_{\ell-1} U$\;
                }
            }
        }
    }
}
\Return{$C^{\ell}_k$, 
$k = 1,\ldots,d, \ \ell = 
1,\ldots,\alpha$.}
\caption{ALS optimization with gauge mediated weight sharing}
\label{alg:als-gmws}
\end{algorithm}

\section{Numerical experiments}\label{sec:numerics}

We will demonstrate our method on three example dynamical systems: The \emph{Fermi-Pasta-Ulam-Tsingou} (FPUT) system, one-dimensional chain of rotating magnetic dipoles and a chain of atoms interacting via a modified Lennard-Jones interaction.

\paragraph*{FPUT system.}
The FPUT system is a chain of non-linear springs with spring constants $\kappa_\ell$.
The dynamical laws are given by
\begin{align}\label{eq:fput}
\begin{split}
f_k(x) =& \kappa_{k+1}(x_{k+1}- x_k) - \kappa_{k} (x_k - x_{k-1})\\
& + \beta_{k+1} (x_{k+1}-x_k)^3 - \beta_{k} (x_k-x_{k-1})^3\\
&k =1,\ldots,d,
\end{split}
\end{align}
which is a one-dimensional system with interaction length $L=1$ and, using the monomial dictionary eq.~\eqref{eq:monomial dict}, separation rank $N=4$, such that each $f_k$ can be represented by TTs with rank bounded by $r=4$.
Furthermore, the polynomial degree of the equations is bounded by $3$ and, hence, we can use TTs \eqref{eq:degree limited TT}, with block-sparse structure given by eq.~\eqref{eq:monomial block sparsity} with $\lambda=3$ and block sizes bounded by $4$, to represent the system exactly.

\paragraph*{Rotating magnetic dipoles.}
Here, we have a chain of magnetic dipoles at positions $X_\ell$ with magnetic dipole moments $M_\ell$ and moments of inertia $I_\ell$.
They are free to rotate in the plane perpendicular to the chain and their angles of rotation are $x_\ell\in[0,2\pi)$.
The dynamical laws are given by 
\begin{equation}\label{eq:magneticdipole}
    f_k(x) = I_kM_k \sum_{\ell\neq k} \frac{M_{\ell}}{|X_k-X_\ell|^3} \sin(x_k-x_\ell).
\end{equation}
The positions are chosen so that $X_1 < X_2 < \dots < X_d$.
In fact, we set $M_\ell = I_\ell = 1, \ X_\ell = \ell - 1 \ \forall \ell \in [d]$.
This is a $2$-mode interacting system, \cref{def:k-mode interacting} with separation rank $N=2$, using the trigonometric dictionary \eqref{eq:trig dict}.
Hence, it suffices to use rank $r=2(d-1)$ TTs.
Since the degree given by $w(1) = 0, w(1) = w(2) = 1$ is bounded by $2$, we can use the block-sparse structure
\begin{align}
\left(C_\ell)\right)_1 = \begin{pmatrix} *&0&0\\0&*&0\\0&0&* \end{pmatrix},&\,\left(C_\ell\right)_2 = \begin{pmatrix} 0&*&0\\0&0&*\\0&0&0 \end{pmatrix},\nonumber\\
\left(C_\ell\right)_3 =& \begin{pmatrix}0&*&0\\0&0&*\\0&0&0 \end{pmatrix}
\end{align}
with block sizes bounded by $2(d-1)$ and represent the laws exactly.
However, this is also a system with $\left(3,\frac1{2\sqrt2\pi}\right)$-algebraically decaying interactions (after rescaling $x_\ell \mapsto (2\pi)^{-1} x_\ell$, so that $x_\ell \in [0,1)$) and separation rank $2$, \cref{def:decaying interactions}.
Hence, \cref{cor:approx locality} allows us to limit the block sizes to a constant (in $d$) and get an approximation of the dynamical law with block-sparse TTs.

\paragraph*{Lennard-Jones chain.}
The final example is a chain of particles of masses $m_\ell$ that interact via a modified Lennard-Jones potential.
This is the hardest example to learn of the three.
The dynamical laws for the positions $x_\ell$ of the particles along the chain are given by 
\begin{align}
    \begin{split}\label{eq:lennardjones}
        f_k(x) =& 6m_k \sum_{\ell \neq k}\sign(x_k-x_\ell) \frac{\varepsilon_{k,\ell}}{R_{{k,\ell}}}\\&\left(2\left(\frac{R_{{k,\ell}}}{|x_k-x_\ell|}\right)^{2q+1}-\left(\frac{R_{{k,\ell}}}{|x_k-x_\ell|}\right)^{q+1}\right), 
    \end{split}
\end{align}
where $\varepsilon_{k,\ell},\ R_{k,\ell}$ are parameters of the interaction between the modes $k,\ell$ and we set $q=2$.
We set $m_k = \varepsilon_{k,\ell} = R_{k,\ell} = 1 \ \forall k,\ell \in [d]$.
Since it is hard to approximate inverse functions with polynomial dictionaries, we learn
\begin{equation}
    g_k(x) \coloneqq (x_k - x_{k-1})^{2q+1} (x_k - x_{k+1})^{2q+1} f_k(x),
\end{equation}
instead of eq.~\eqref{eq:lennardjones} from the accordingly transformed data $\{x_i,\tilde y_i\}_{i\in[M]}$.

This is a $2$-mode interacting system with separation rank with respect to polynomial dictionaries $N = p^2$, since polynomial expansions of inverse functions contain an infinite number of terms.
By \cref{thrm:TT rep of K-mode sys}, we can represent the dynamical laws by TTs with rank bounded by $r=p^2 (d-1)$.
For good approximations of the inverse function, we require large $p$, so for practical use-cases we would like to limit the rank more.
This is justified by \cref{cor:approx locality}, since this is also a system with $(2,g)$-algebraically decaying interactions and separation rank $p^2$, if there is a finite amount of energy in the system and the initial conditions are chosen so that we can ensure that the state space is $\mathcal{S} \subset \bigtimes_{\ell\in[d]}\mathcal{I}_\ell$, where $\mathcal{I}_\ell$ are finite intervals for all $\ell \in [d]$.

The polynomial degree of eq.~\eqref{eq:lennardjones} is bounded by $2p$, such that we can use block-sparsity, where each core $\left(C_\ell\right)_i$ has non-zero blocks only on the $i$-th diagonal.

For completeness, we include the formula for the total energy in the system
\begin{equation}
\begin{split}
    E = \sum_{k < \ell}& \varepsilon_{{k,\ell}}\Big[\Big(\frac{R_{{k,\ell}}}{|x_k(0)-x_\ell(0)|}\Big)^{2q}-\\
    &\Big(\frac{R_{{k,\ell}}}{|x_k(0)-x_\ell(0)|}\Big)^q\Big] +\sum_k \frac{1}{2} m_k\dot{x}_k^2(0).
\end{split}
\end{equation}

\subsection{Results}\label{ssec:results}

For each of the three example systems, we randomly draw data $\mathcal{D} = \{x_i\}_{i\in[M]}$ from the corresponding $\mathcal{S}$ and compute $y_i = f(x_i)$ or $y_i = f(x_i) + \eta_i^\sigma$, where the elements of $\eta_i^\sigma$ are drawn from the Gaussian distribution with standard deviation $\sigma$.
Unless stated otherwise, we use noiseless data with $\sigma=0$.
Therefore, the data that we use for learning does not come from trajectories of the dynamical systems, but instead they are random (potentially noisy) evaluations of the dynamical law $f(x)$.
This somewhat simplifies the setting, especially since we do not have to approximate $\dot x(t)$ or $\ddot x(t)$ in eq.~\eqref{eq:dynamical law} from the trajectory, using e.g. finite differences.
However, since we can hope to learn the dynamical laws only from trajectories that sufficiently explore the state space $\mathcal{S}$, sampling $\mathcal{S}$ at random is not too different from using such trajectories.

Given $\mathcal{D}$, we use ALS-GMWS to find a $S_L$-self-similar low rank block-sparse TT representation of an estimate $\hat f$ of the dynamical law.
We benchmark the quality of the estimate with respect to the true dynamical law $f$ via the residuum
\begin{equation}
    \operatorname{res}(\hat f,f) \coloneqq \sqrt{\frac{\sum_{i\in[M']}\left\|\hat f(x'_i) - f(x'_i)\right\|_2^2}{\sum_{j\in[M']}\|f(x'_j)\|_2^2}},
\end{equation}
where $\{x'_i\}_{i\in[M']}$ are random samples from $\mathcal{S}$, which are different to the samples used for training.
In particular, we use $M' = 2\times10^4$.

All experiments have been conducted on consumer grade hardware and the code has not been optimized for speed.

The block-sparsity used for polynomial dictionaries is that non-zero blocks of $\left(C_\ell\right)_i$ are on the $i$-th block-diagonal, while for the trigonometric dictionary eq.~\eqref{eq:trig dict} it is such that $\left(C_\ell\right)_0$ is block-diagonal, while $\left(C_\ell\right)_1$ and $\left(C_\ell\right)_2$ have non-zero blocks on the first block-diagonal.
In both cases, if we bound the degree by $\lambda$, the number of block-rows and block-columns is $\lambda+1$.
The maximum block size is $\rho$.
For all experiments we have used self-similarity given by the selection table $S_L$ defined in eq.~\eqref{eq:local selection table}.

\begin{figure}
\centering
\includegraphics{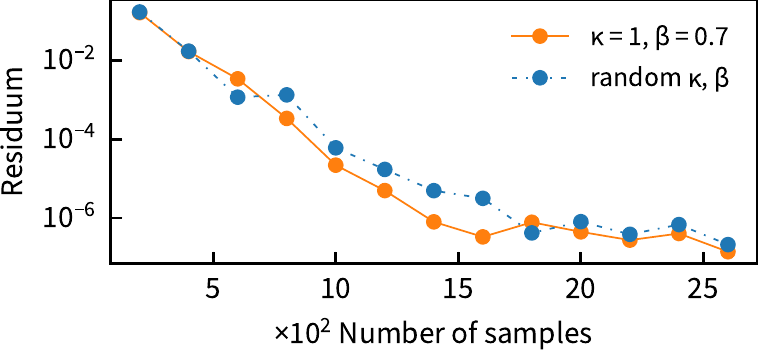}
\caption{\textbf{FPUT system.} Recovery of a $d = 50$ system with constant spring constants (blue crosses) and spring constants $\kappa, \ \beta$ drawn from the uniform distribution on $[0,2]$ and $[0,1.4]$ respectively (orange plus signs). 
Residuum after $10$ ALS-GMWS sweeps is plotted against the size of the learning set.
We use degree $3$ Legendre polynomial dictionary, maximum block size $\rho = 2$ and selection table interaction length $L=5$.}
\label{fig:fput}
\end{figure}

\paragraph*{FPUT system.}
To recover FPUT systems eq.~\eqref{eq:fput} with size $d=50$, we use degree $3$ Legendre polynomial dictionary.
In \cref{fig:fput} we show the recovery of the translationally invariant FPUT system and FPUT system with randomly sampled spring constants.
We plot the residuum achieved after $10$ ALS-GMWS sweeps, using varying numbers $M$ of training samples.
Both systems are successfully recovered using around $2\times 10^3$ training samples.

\paragraph*{Rotating magnetic dipoles.}
We perform three experiments on the rotating magnetic dipole chain.

\begin{figure}
\centering
\includegraphics{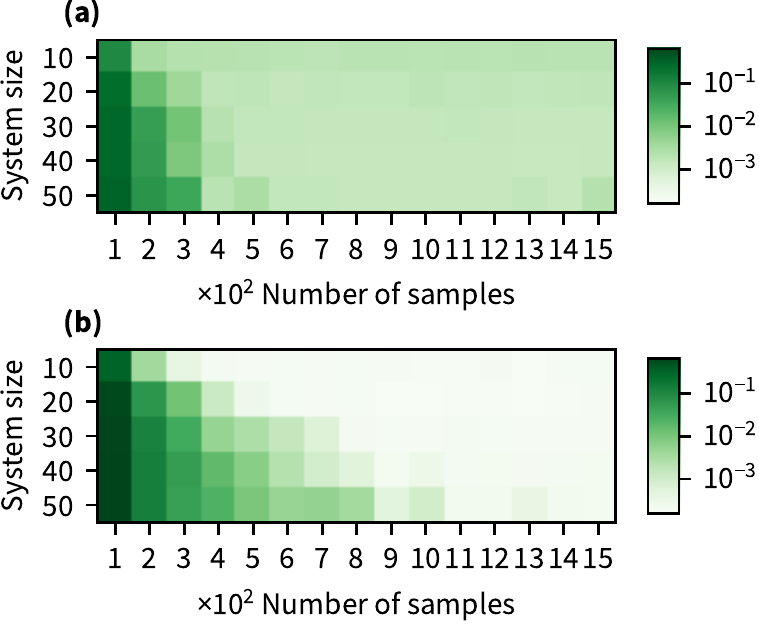}
\caption{\textbf{Magnetic dipole chain with trigonometric dictionary.} Residuum after $8$ ALS-GMWS sweeps on varying system sizes $d$ as a function of the size of the learning set.
The maximum block size is set to $\rho=3$ and the selection table interaction length to \textbf{(a)} $L=5$ and \textbf{(b)} $L=9$.
}
\label{fig:magnetic_trig}
\end{figure}

First, we use the trigonometric dictionary to recover chains with $10 \le d \le 50$.
In \cref{fig:magnetic_trig} we plot the residua achieved using varying numbers $M$ of training samples.
We show the results for $L=5,9$.
For $d=50$, successful recovery requires around $4\times 10^2$ training samples for $L=5$ and $9\times10^2$ training samples for $L=9$, which however achieves around $10$ times smaller residuum.

\begin{figure}
\centering
\includegraphics{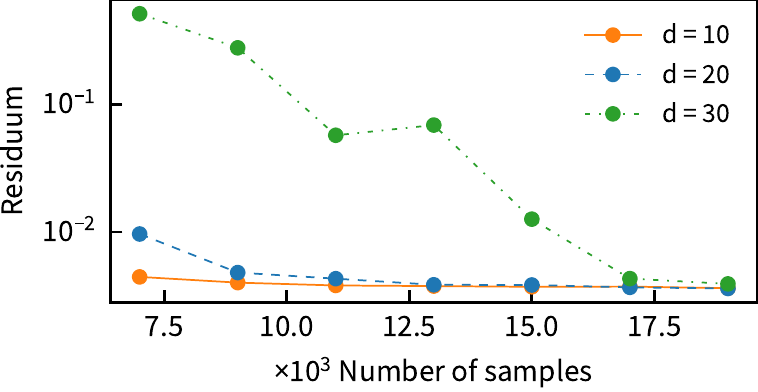}
\caption{\textbf{Magnetic dipole chain with Legendre polynomial dictionary.}
Residuum after $8$ ALS-GMWS sweeps on system sizes $d=10,20,30$ is plotted as a function of the learning set size.
The maximum block size is set to $\rho=3$, dictionary is truncated at degree $9$ and selection table interaction length to $L=5$.}
\label{fig:magnetic_legendre}
\end{figure}

Second, we use degree $9$ Legendre polynomial dictionary to recover chains with $d=10,20,30$.
The residua for varying numbers $M$ of training samples are shown in \cref{fig:magnetic_legendre}.
For $d=30$, we require around $1.7 \times 10^4$ samples.
This demonstrates the importance of choosing an appropriate dictionary, when compared with the previous experiment.

\begin{figure}
\centering
\includegraphics{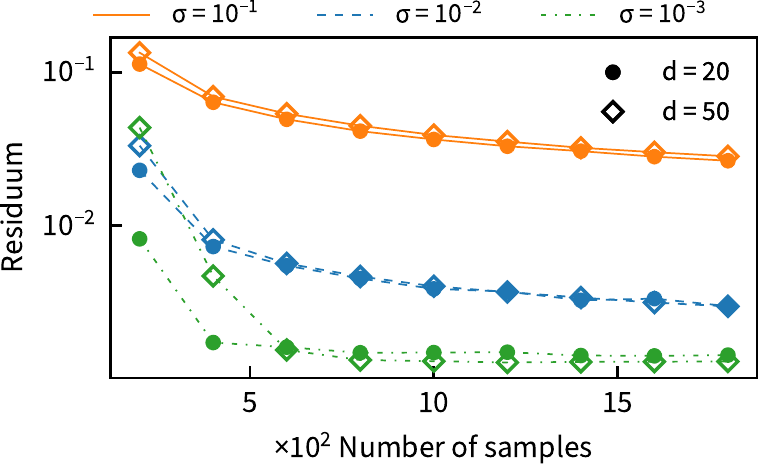}
\caption{\textbf{Magnetic dipole chain from noisy data.}
Residuum after $8$ ALS-GMWS sweeps using noisy data with varying levels of noise $\sigma$ for system sizes $d=20$ (crosses) and $d=50$ (plus signs). 
We use the trigonometric dictionary, maximum block size $\rho=3$ and selection table interaction length $L=5$.}
\label{fig:magnetic_noise}
\end{figure}

Finally, we use trigonometric dictionary to recover chains with $d=20,50$ from noisy data with varying levels $\sigma$ of noise.
We plot the residua achieved using varying numbers $M$ of training samples in \cref{fig:magnetic_noise}.
The results show recovery of the system down to the noise level, demonstrating noise robustness of the proposed method.

\begin{figure}
\centering
\includegraphics{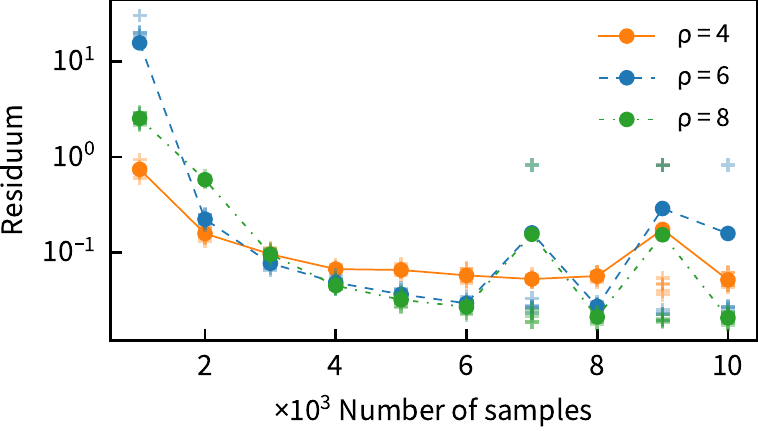}
\caption{\textbf{Lennard-Jones chain.}
The residuum after $8$ sweeps of the ALS-GMWS algorithm is plotted for system size $d=10$ and maximum block sizes $\rho=4,6,8$.
The mean (crosses) over 6 runs of the algorithm is plotted, as well as the result of the individual instances (plus signs).
We use degree $8$ Legendre polynomial dictionary and set selection table interaction length to $L = 5$.
}
\label{fig:lennard_jones}
\end{figure}

\paragraph*{Lennard-Jones chain.}
We recover the Lennard-Jones chain eq.~\eqref{eq:lennardjones} for $d=10$, using degree $8$ Legendre polynomial dictionary.
We set $L=5$.
\cref{fig:lennard_jones} shows the residua after $8$ ALS-GMWS sweeps for various maximum block sizes, as a function of the training set size.
This is clearly the hardest example, requiring around $6\times 10^3$ samples for successful recovery of even such a small system.
Furthermore, the plot shows three initialization instances when the algorithm has not converged well.
In practice it is therefore sometimes beneficial to run the algorithm multiple times with different initializations.

\section{Conclusions}\label{sec:conclusions}

Learning dynamical laws is a key task that has been moving into the focus of attention. 
The well-known and immensely popular SINDy approach introduces a data-driven algorithm for obtaining dynamical systems from data.
This approach allows for a reliable recovery for small systems, but is not scalable to a large number of degrees of freedom: For this, a meaningful, physically motivated restriction of the hypothesis class is necessary.
In this work, we have overcome this obstacle.
We have shown that block-sparse \emph{tensor trains} (TT) with self-similarity provide a suitable efficient ansatz class for learning dynamical laws from data in many contexts of practical interest.
In particular, we have proven that these include local one-dimensional systems, one-dimensional systems with algebraically decaying interactions and systems with $K$-body interactions in any number of dimensions.
For learning dynamical laws within this class, we have developed a new variant of the \emph{alternating least squares} (ALS) algorithm for block-sparse TTs, which we refer to as \emph{ALS with gauge mediated weight sharing} (ALS-GMWS), which is suitable for self-similar systems.
The method has been successfully demonstrated on three physically relevant one-dimensional dynamical systems and robustness to Gaussian additive noise in the data has been demonstrated.

\section{Acknowledgements}
We thank A.~Goe{\ss}mann and R.~Schneider for insightful discussions. This work has been funded by the DFG (Interdisciplinary Research Training Group DAEDALUS RTG 2433 and EI 519/15-1), Ei 519/20, and the Cluster of Excellence MATH+. it has also been funded by the ERC (DebuQC).

\section{Code availability}
The code used in this publication, written in the Python language, is available in the repository \url{https://github.com/JonasFuksa/bstt}. The code has been developed by M.~Götte and P.~Trunschke. J.~Fuksa performed the numerical experiments.


\bibliographystyle{alpha}
\newcommand{\etalchar}[1]{$^{#1}$}

\onecolumn 

\appendix

\section{Proof of Theorem \ref{thrm:approximate locality}}\label{app:approximate locality proof}

\cref{thrm:approximate locality} bounds for $\chi>1$ the error that we obtain by representing the dynamical law of a system with $(\chi,g)$-algebraically decaying interactions and separation rank $N$ by a low rank TT.
The original system is governed by the dynamical law
\begin{equation}
    f_k(x) = \sum_{L \in[d]} L^{-\chi} g_{k,L}(x).
\end{equation}
Suppose the approximation of this system
\begin{equation}\label{eq:truncated algebraic sys}
    \tilde f_k^{(\tilde L)}(x) = \sum_{L \in [\tilde L]} L^{-\chi} g_{k,L}(x).
\end{equation}
The approximation error is bounded by
\begin{align}
    \|f_k(x) - \tilde f_k^{(\tilde L)}(x)\|_2 &\le \sum_{L = \tilde L + 1}^d L^{-\chi} g\\
    &\le \sum_{L = \tilde L + 1}^\infty L^{-\chi} g\\
    &\le c_1(\chi,\tilde L) g\\
    &= \mathcal{O}(L^{1-\chi}) g,
\end{align}
where we use $L_2$ norm sub-additivity to get the first inequality and the last inequality comes from using the upper bound $\sum_{j=N}^\infty f(j) \le f(N) + \int_{N}^\infty f(x) \mathrm{d}x$ for $f(x)$ that decreases monotonously for all $x\in[N,\infty)$.
Evaluating the integral for the present case $N = \tilde L + 1$, $f(x) = x^{-\chi}$ with $\chi>1$, we get
\begin{align}
    \sum_{L=\tilde L+1}^\infty L^{-\chi} &\le (\tilde L + 1)^{-\chi} \frac{\tilde L + \chi}{\chi-1}\\
    &\eqqcolon c_1(\chi,\tilde L)\\
    &= \mathcal{O}(\tilde L^{1-\chi}).
\end{align}

It remains to be shown that the truncated dynamical system corresponding to the dynamical law $\tilde f^{(\tilde L)}$ is one-dimensional with interaction length $\tilde L$ and separation rank $N\tilde L$. 
Since all of the terms in the sum eq.~\eqref{eq:truncated algebraic sys} depend trivially on $x_i$ with $i \notin [k-\tilde L, k+\tilde L]$, the system has interaction length $\tilde L$ by definition.
Moreover, since it is a sum of $\tilde L$ functions with separation rank at most $N$ with respect to any bipartition $\mathcal{P}_k$, the total separation rank is bounded by $N \tilde L$.
$\square$

\section{Proof of Corollary \ref{cor:approx locality}}\label{app:approx locality corollary proof}

\cref{cor:approx locality} shows that there exists $\varepsilon$-approximate low rank TT representation of a dynamical systems with $(\chi,g)$-algebraically decaying interactions and separation rank $N$, if $\chi>1$.
To prove this, use \cref{thrm:approximate locality} to find that there exists a one-dimensional dynamical system with interaction length $\tilde L$ and separation rank $N\tilde L$, which is governed by a dynamical law $\tilde f$, such that, for all $k\in[d]$, $\tilde f_k$ is an $\varepsilon$-approximation of $f_k$ in the $L_2$ norm, if
\begin{align}
    (\tilde L + 1)^{-\chi} \frac{\tilde L + \chi}{\chi-1} &\le \frac\varepsilon g\\
    (\tilde L + 1)^{1-\chi} + (\chi-1)(\tilde L + 1)^{-\chi} &\le (\chi-1)\frac\varepsilon g.
\end{align}
This will certainly be satisfied if 
\begin{equation}
    (\tilde L + 1)^{1-\chi} \le \frac{\chi-1}{\chi}\frac\varepsilon g.
\end{equation}
Rearranging, we obtain
\begin{equation}
    \tilde L \ge \left(\frac\chi{\chi-1} \frac g\varepsilon\right)^{\frac1{\chi-1}} - 1.
\end{equation}
From \cref{thrm:TT rep of local systems} we know that in this approximation, each $\tilde f_k$ admits a TT representation with rank $r \ge N \tilde L$, so 
\begin{equation}
    r \ge N \left[\left(\frac{\chi-1}\chi \frac\varepsilon g\right)^{\frac1{1-\chi}} - 1\right].
\end{equation}
$\square$

\section{Proof of Theorem \ref{thrm:TT rep of K-mode sys}}\label{app:K-mode proof}

\cref{thrm:TT rep of K-mode sys} bounds the TT rank of the TT representation of the dynamical law $f(x)$. 
In order to do this, we need to bound the separation rank with respect to bipartitions $\mathcal{P}_k = (\{x_1,\dots,x_k\},\{x_{k+1},\dots,k_d\})$. 
To prove 1., first we need to show that 
\begin{align}
\mathcal{K}_k \le \genfrac(){0pt}{0}{d-1}{K-1}.
\end{align}
To see this note that we require each $\mathcal{J}_\ell^{(k)}$ to be distinct and to contain $k$.
Hence, the upper bound is the number of ways of selecting the remaining $K-1$ elements of $\mathcal{J}_\ell^{(k)}$.
Now, rewrite the decomposition of the dynamical law of a $K$-mode interacting system with separation rank $N$ as
\begin{align}
    f_k(x) &= \sum_{\ell \in [\mathcal{K}_k]} \sum_{(i_{j_1},\dots,i_{j_K}) \in \mathcal{I}_{\ell}^{(k)}} g_{i_{j_1}}(x_{j_1}) \dots g_{i_{j_K}}(x_{j_K})\\
    &= \sum_{\ell \in [\mathcal{K}_k]}\sum_{(i_{j_1},\dots,i_{j_K})\in \mathcal{I}_{\ell}^{(k)}} \bigotimes_{m=1}^d \begin{cases}
        g_{i_{m}}(x_{m}) & m\in\mathcal{J}_\ell^{(k)},\\
        \id(x_{j_m}) & \text{otherwise},
    \end{cases}
    \nonumber
\end{align}
where $\id(x) = 1$.
This is a sum of tensor products of single variable functions of each mode.
Hence, the separation rank of $f_k$ with respect to any $\mathcal{P}_k$ is upper bounded by the number of terms in the sum, which is $\mathcal{K}_k|\mathcal{I}_{k}^{(\ell)}| \le \genfrac(){0pt}{0}{d-1}{K-1} N$.
Since the minimal ranks of the TT representation of $f_k$ are equal to the corresponding separation ranks, we get that 
\begin{equation}
    r_k \le r = N \genfrac(){0pt}{0}{d-1}{K-1} = N\mathcal{O}(d^{K-1}).
\end{equation}

The proof of 2. proceeds along similar lines. 
We write the total dynamical law $f(x)$, which we think of in terms of eq.~\eqref{eq:total TT}, as
\begin{equation}
    f(x) = \sum_{k' \in [d]} \sum_{\ell \in [\mathcal{K}_{k'}]} \tilde g_{k',\ell}(x) \otimes e_{k'},
\end{equation}
where $e_{k'} \in \mathbb{R}^d$ is the vector with one at the $k'$-th element and zeroes elsewhere.
Here each $\tilde g_{k',\ell}$ represents a function that depends non-trivially only on $x_i$ with $i \in \mathcal{J}_\ell^{(k')}$ and which has a separation rank at most $N$ with respect to any bipartition of $[d]$, due to the assumption that $f(x)$ is a dynamical law of a $K$-body interacting system with separation rank $N$.
Given a bipartition $\mathcal{P}_k = (\mathcal{P}^\text{left},\mathcal{P}^\text{right}) = \{\{x_1,\dots,x_k\},\{x_{k+1},\dots,x_d\}\}$, we can write
\begin{align}
    f(x) &= \id_{1,k} \otimes \sum_{k' \in [d]} \sum_{\ell \in [\mathcal{K}_{k'}]: \mathcal{J}_\ell^{(k')} \subset \mathcal{P}_k^\text{right}} \tilde g'_{k',\ell} \otimes e_{k'} + \label{eq:term one} \\
    &\quad + \sum_{k'\in[d]} \sum_{\ell \in [\mathcal{K}_{k'}]: \mathcal{J}_\ell^{(k')} \subset \mathcal{P}_k^\text{left}} \tilde g''_{k',\ell} \otimes \id_{k+1,d} \otimes e_{k'} + \label{eq:term two}  \\
    &\quad + \sum_{k'\in[d]} \sum_{\ell \in [\mathcal{K}_{k'}]: \mathcal{J}_\ell^{(k')} \nsubseteq \mathcal{P}_k^\text{left,right}} \tilde g_{k',\ell} \otimes e_{k'},\label{eq:term three}
\end{align}
where the $\id_{k_1,k_2}$ is the one function of $x_{k_1},\dots,x_{k_2}$ and we abuse the subset notation in $\mathcal{J}_\ell^{(k')} \subset \mathcal{P}_k^{\text{left,right}}$ and similar to indicate that $x_j \in \mathcal{P}_k^{\text{left,right}}$ for all $j \in \mathcal{J}_\ell^{(k')}$.
Furthermore, we denote by $g'_{k',\ell}$, $g''_{k',\ell}$ the restriction of $g_{k',\ell}$ onto the modes $x_{k+1},\dots,x_d$ and $x_1,\dots,x_k$ respectively.
Note that this is well defined since we always use this notation in the cases where $g$ depends trivially on the modes that we throw away.

Each term \eqref{eq:term one}, \eqref{eq:term two}, \eqref{eq:term three} is now written in such a way that we can read of a bound on its separation rank with respect to $\mathcal{P}_k$.
The term \eqref{eq:term one} has separation rank with respect to $\mathcal{P}_k$ at most $1$.
The term \eqref{eq:term two} vanishes if $k' \ge k+1$, since then the condition on the second sum cannot be satisfied. 
For each $k'\le k$, the number of terms in the sum, each of which has separation rank $1$ with respect to $\mathcal{P}_k$ is upper bounded by $\genfrac(){0pt}{1}{k-1}{K-1}$, so the bound on the separation rank of eq.~\eqref{eq:term two} is $k\genfrac(){0pt}{1}{k-1}{K-1}$.
Finally, in the last term \eqref{eq:term three}, we know that each $\tilde g_{k',\ell}$ has separation rank with respect to any bipartition bounded by $N$.
The number of $\ell \in [\mathcal{K}_{k'}]$ such that $\mathcal{J}_\ell^{(k')}$ satisfies the condition is
\begin{equation}
    \begin{split}
        c_2(N,d,k) &\coloneqq N\left[d \genfrac(){0pt}{0}{d-1}{K-1} - k \genfrac(){0pt}{0}{k-1}{K-1} - (d-k) \genfrac(){0pt}{0}{d-k-1}{K-1}\right] \\
        &= \mathcal{O}(Nd^K).
    \end{split}
\end{equation}
Putting all the bounds together, we recover claim 2. of the theorem.
$\square$

\section{Proof of Theorem \ref{thrm:block sparsity}}\label{app:block sparsity proof}

\cref{thrm:block sparsity} states that functions with fixed degree given by the degree map $w$ admit a TT representation, such that the left and right interface tensors satisfy the eigenvalue equations
\begin{align}
    \phi^{>\ell}L^{>\ell} &= \Lambda^{>\ell}\phi^{>\ell}\label{eq:block sparsity 1 in proof},\\
    L^{<\ell+1}\phi^{<\ell+1} &= \phi^{<\ell+1}\left(\lambda \id - \Lambda^{>\ell}\right)\label{eq:block sparsity 2 in proof}.
\end{align}
We have a TT $\phi$ with tensor cores $\{C_i\}_{i\in[d]}$, such that it can be written in tensor network notation as 
\begin{equation}
    \begin{tikzpicture}[
        dot/.style = {circle,fill,minimum size=#1,inner sep=0pt,outer sep=0pt},
        dot/.default = 5pt
        ]
        \node[dot,label=above:$C_1$] at (0,0) (c1) {};
        \node[dot,label=above:$C_2$] [right=.5cm of c1] (c2) {};
        \node [right=.5cm of c2] (dots) {$\dots$};
        \node[dot,label=above:$C_d$] [right=.5cm of dots] (cd) {};
        \node [left=.5 of c1] (lhs) {$\phi = $};

        \node [below=.3cm of c1] (d1) {};
        \node [below=.3cm of c2] (d2) {};
        \node [below=.3cm of cd] (dd) {};

        \draw (d1) -- (c1) -- (dots.west);
        \draw (d2) -- (c2) -- (dots.west);
        \draw (dd) -- (cd) -- (dots.east);
        \node [right=.2 of cd] (dot) {.};
    \end{tikzpicture}
\end{equation}
Without loss of generality, we can assume that this TT is in left-canonical form, so that the tensor cores satisfy eq.~\eqref{eq:left canonical}, and that the ranks are minimal.
If this is not the case, we can always find a gauge transformation that puts the TT into this form.
We also assume that $\phi$ is an eigenvector of $L$, such that
\begin{equation}\label{eq:eval bs}
    L\phi = \lambda \phi.
\end{equation}
To prove the theorem, we will inductively gauge transform each tensor core, starting at $C_d$ and proceeding one-by-one towards $C_1$, such that after transforming $C_\ell$, for all $\ell'\ge\ell-1$ the right interface tensor satisfies
\begin{equation}
    \phi^{>\ell'}L^{>\ell'} = \Lambda^{>\ell'}\phi^{>\ell'}.
\end{equation}
Finally, we will show that eq.~\eqref{eq:block sparsity 2 in proof} follows from eq.~\eqref{eq:block sparsity 1 in proof}.

\emph{Base case.} First, we will find an appropriate gauge transformation for $C_d$. 
We can write the eigenvalue equation in the form
\begin{equation}
    \begin{tikzpicture}[
        dot/.style = {circle,fill,minimum size=#1,inner sep=0pt,outer sep=0pt},
        dot/.default = 5pt
        ]
        \node (lhs) at (0,0) {$L \phi = \lambda \phi = $};
        \node [dot,label=above:$\phi^{<d}$] (phi) [right=.5 of lhs] {};
        \node [dot,label=above:$C_d$] (cd) [right=.5 of phi] {};
        \node [dot,label=left:$L^{<d}$] (Ll) [below=.3 of phi] {};
        \node (Llb) [below=.3 of Ll]{};
        \node (id1) [below=.8 of cd] {};

        \draw[double] (Llb) -- (Ll) -- (phi);
        \draw (phi) -- (cd) -- (id1);
        
        \node (plus) [right=.5 of cd] {$+$};

        \node [dot,label=above:$\phi^{<d}$] (phi2) [right=.5 of plus] {};
        \node [dot,label=above:$C_d$] (cd2) [right=.5 of phi2] {};
        \node [dot,label=right:$\Omega$] (Lr) [below=.3 of cd2] {};
        \node (id2) [below=.8 of phi2] {};
        \node (Lrb) [below=.3 of Lr] {};
        \node (dot) [right=.3 of cd2] {,};

        \draw[double] (id2) -- (phi2);
        \draw (phi2) -- (cd2) -- (Lr) -- (Lrb);
    \end{tikzpicture}
\end{equation}
where the double line combines multiple indices into a single edge and $\Omega = \diag(w(1),\dots,w(d))$.
We can now contract the first $d-1$ physical indices with $\left(\phi^{<d}\right)^*$ and use the assumption that the TT is written in left-canonical form to obtain
\begin{equation}
    \begin{tikzpicture}[
        dot/.style = {circle,fill,minimum size=#1,inner sep=0pt,outer sep=0pt},
        dot/.default = 5pt
        ]
        \node (lhs) at (0,0) {$\lambda$};
        \node [dot,label=above:$C_d$] (cd1) [right=.5 of lhs] {};
        \node (cdb) [below=.5 of cd1] {};

        \draw (lhs.east) -- (cd1) -- (cdb);

        \node (eq) [right=.3 of cd1] {$=$};

        \node [dot,label=above:$\phi^{<d}$] (phi) [right=.5 of eq] {};
        \node [dot,label=above:$C_d$] (cd) [right=.5 of phi] {};
        \node [dot,label=left:$L^{<d}$] (Ll) [below=.3 of phi] {};
        \node (id1) [below=.8 of cd] {};
        \node [dot,label=left:$\left(\phi^{<d}\right)^*$] (phi*) [below=.5 of Ll] {};
        \node (right) [right=.3 of phi*] {};

        \draw[double] (phi*) -- (Ll) -- (phi);
        \draw (phi) -- (cd) -- (id1);
        \draw (phi*) -- (right);
        
        \node (plus) [right=.5 of cd] {$+$};

        \node [dot,label=above:$C_d$] (cd2) [right=.5 of plus] {};
        \node [dot,label=right:$\Omega$] (Lr) [below=.3 of cd2] {};
        \node (left) [left=.3 of cd2] {};
        \node (Lrb) [below=.3 of Lr] {};
        \node (dot) [right=.3 of cd2] {,};

        \draw (left) -- (cd2) -- (Lr) -- (Lrb);
    \end{tikzpicture}
\end{equation}
which, rearranging, we can write in matrix notation as 
\begin{equation}\label{eq:base case eval eq}
    \left(\phi^{<d}\right)^\dagger L^{<d} \phi^{<d} C_d = C_d(\lambda\id - \Omega).
\end{equation}
Since $\left(\phi^{<d}\right)^\dagger L^{<d} \phi^{<d}$ is Hermitian, there exists a unitary $U_{d}$ and a diagonal matrix $\Lambda^{<d}$ with non-decreasing diagonal entries, such that $\left(\phi^{<d}\right)^\dagger L^{<d} \phi^{<d} = U_{d}^\dagger \Lambda^{<d} U_{d}$.
Now we can write eq.~\eqref{eq:base case eval eq} as 
\begin{equation}\label{eq:base case block sparsity}
    \Lambda^{<d} \tilde C_d = \tilde C_d (\lambda \id - \Omega),
\end{equation}
where $\tilde C_d = U_d C_d$.
Since $U_d$ is a unitary, it defines a gauge transformation 
\begin{equation}
    C_{d-1} \mapsto \hat C_{d-1} = C_{d-1} U_{d-1}^\dagger \ , \quad C_d \mapsto \tilde C_d = U_{d-1} C_d,
\end{equation}
which leaves $\phi$ invariant and preserves its left-canonical form.
Rewriting eq.~\eqref{eq:base case block sparsity} as $\tilde C_d \Omega = (\lambda\id - \Lambda^{<d})\tilde C_d$ and noticing that $\Omega = L^{>d-1}$ and $\tilde C_d = \phi^{>d-1}$, we get
\begin{equation}
    \phi^{>d-1}L^{>d-1} = \Lambda^{>d-1}\phi^{>d-1},
\end{equation}
where we defined $\Lambda^{>d-1} = \lambda\id - \Lambda^{<d}$, which has non-increasing entries.
This is eq.~\eqref{eq:block sparsity 1 in proof} for $\ell=d-1$.

\emph{Induction step.} Take $\ell \in \{2,\dots,d-1\}$ and assume that $\phi$ is in left canonical form with 
\begin{equation}\label{eq:inductive hypothesis}
    \phi^{>\ell'}L^{>\ell'} = \Lambda^{>\ell'}\phi^{>\ell'} \quad \forall \ell'\ge\ell,
\end{equation}
where $\Lambda^{>\ell}$ is a diagonal matrix with non-increasing entries.
We can decompose the LHS of eq.~\eqref{eq:eval bs} to get 
\begin{equation}
    \begin{tikzpicture}[
        dot/.style = {circle,fill,minimum size=#1,inner sep=0pt,outer sep=0pt},
        dot/.default = 5pt
        ]
        \node (lhs) at (0,0) {$\lambda \phi = L \phi =$};
        \node [dot,label=above:$\phi^{<\ell}$] (phi) [right=.5 of lhs] {};
        \node [dot,label=above:$C_\ell$] (cl) [right=.5 of phi] {};
        \node [dot,label=left:$L^{<\ell}$] (Ll) [below=.3 of phi] {};
        \node (Llb) [below=.3 of Ll]{};
        \node (id1) [below=.8 of cl] {};
        \node [dot,label=above:$\phi^{>\ell}$] (phiR) [right=.5 of cl] {};
        \node (idr1) [below=.8 of phiR] {};

        \draw[double] (Llb) -- (Ll) -- (phi);
        \draw (phi) -- (cl) -- (id1);
        \draw (cl) -- (phiR);
        \draw[double] (phiR) -- (idr1);
        
        \node (plus) [right=.5 of phiR] {$+$};

        \node [dot,label=above:$\phi^{<\ell}$] (phi2) [right=.5 of plus] {};
        \node [dot,label=above:$C_\ell$] (cl2) [right=.5 of phi2] {};
        \node [dot,label=left:$\Omega$] (omega) [below=.3 of cl2] {};
        \node (Llb2) [below=.3 of omega]{};
        \node (idl2) [below=.8 of phi2] {};
        \node [dot,label=above:$\phi^{>\ell}$] (phiR2) [right=.5 of cl2] {};
        \node (idr2) [below=.8 of phiR2] {};

        \draw[double] (idl2) -- (phi2);
        \draw (phi2) -- (cl2) -- (omega) -- (Llb2);
        \draw (cl2) -- (phiR2);
        \draw[double] (phiR2) -- (idr2);

        \node (plus2) [right=.5 of phiR2] {$+$};

        \node [dot,label=above:$\phi^{<\ell}$] (phi2) [right=.5 of plus2] {};
        \node [dot,label=above:$C_\ell$] (cl2) [right=.5 of phi2] {};
        \node (idl2) [below=.8 of phi2] {};
        \node [dot,label=above:$\phi^{>\ell}$] (phiR2) [right=.5 of cl2] {};
        \node [dot,label=right:$L^{>\ell}$] (lr) [below=.3 of phiR2] {};
        \node (Llb2) [below=.3 of lr]{};
        \node (idr2) [below=.8 of cl2] {};
        \node (dot) [right=.5 of phiR2] {.};

        \draw[double] (idl2) -- (phi2);
        \draw (phi2) -- (cl2) -- (idr2);
        \draw (cl2) -- (phiR2);
        \draw[double] (phiR2) -- (lr) -- (Llb2);
    \end{tikzpicture}
\end{equation}
Contracting the first $\ell-1$ physical indices with $\left(\phi^{<\ell}\right)^*$, using the left-canonical gauge condition and eq.~\eqref{eq:inductive hypothesis}, when we fix the $\ell$-th physical index to $i$, we obtain
\begin{equation}
    \begin{tikzpicture}[
        dot/.style = {circle,fill,minimum size=#1,inner sep=0pt,outer sep=0pt},
        dot/.default = 5pt
        ]
        \node (lmbd) at (0,0) {$\lambda$};
        \node [dot,label=above:$C_\ell$] (cl1) [right=.5 of lmbd] {};
        \node (cl1l) [left=.3 of cl1] {};
        \node (i1) [below=.3 of cl1] {$i$};
        \node [dot,label=above:$\phi^{>\ell}$] (phiR1) [right=.5 of cl1] {};
        \node (phiR1b) [below=.5 of phiR1] {};

        \draw (cl1l) -- (cl1) -- (i1);
        \draw (cl1) -- (phiR1);
        \draw[double] (phiR1) -- (phiR1b);

        \node (eq) [right=.5 of phiR1] {$=$};

        \node [dot,label=above:$\phi^{<\ell}$] (phi) [right=.5 of eq] {};
        \node [dot,label=above:$C_\ell$] (cl) [right=.5 of phi] {};
        \node [dot,label=left:$L^{<\ell}$] (Ll) [below=.2 of phi] {};
        \node [dot,label=below:$\left(\phi^{<\ell}\right)^*$] (phi*) [below=.2 of Ll] {};
        \node (phi*r) [right=.2 of phi*] {};
        \node (i2) [below=.3 of cl] {$i$};
        \node [dot,label=above:$\phi^{>\ell}$] (phiR) [right=.5 of cl] {};
        \node (idr1) [below=.5 of phiR] {};

        \draw[double] (phi*) -- (Ll) -- (phi);
        \draw (phi) -- (cl) -- (i2);
        \draw (cl) -- (phiR);
        \draw[double] (phiR) -- (idr1);
        \draw (phi*) -- (phi*r);
        
        \node (plus) [right=.5 of phiR] {$+ \quad w(i)$};

        \node [dot,label=above:$C_\ell$] (cl2) [right=.4 of plus] {};
        \node (cl2l) [left=.3 of cl2] {};
        \node (Llb2) [below=.3 of omega]{};
        \node (i3) [below=.3 of cl2] {$i$};
        \node [dot,label=above:$\phi^{>\ell}$] (phiR2) [right=.5 of cl2] {};
        \node (idr2) [below=.5 of phiR2] {};

        \draw (cl2l) -- (cl2) -- (i3);
        \draw (cl2) -- (phiR2);
        \draw[double] (phiR2) -- (idr2);

        \node (plus2) [right=.5 of phiR2] {$+$};

        \node [dot,label=above:$\phi^{<\ell}$] (phi2) [right=.5 of plus2] {};
        \node [dot,label=above:$C_\ell$] (cl3) [right=.5 of phi2] {};
        \node (idl2) [below=.5 of phi2] {};
        \node [dot,label=above:$\Lambda^{>\ell}$] (Lmbd) [right=.5 of cl3] {};
        \node [dot,label=above:$\phi^{>\ell}$] (phiR2) [right=.5 of Lmbd] {};
        \node (idphiR) [below=.5 of phiR2] {};
        \node (i4) [below=.3 of cl3] {$i$};
        \node (dot) [right=.5 of phiR2] {,};

        \draw[double] (idl2) -- (phi2);
        \draw (phi2) -- (cl3) -- (i4);
        \draw (cl3) -- (phiR2);
        \draw[double] (phiR2) -- (idphiR);
    \end{tikzpicture}
\end{equation}
which we can rearrange and write in matrix notation as 
\begin{equation}\label{eq:induction eval eq}
    \left(\phi^{<\ell}\right)^\dagger L^{<\ell} \phi^{<\ell} \left(C_\ell\right)_i = \left(C_\ell\right)_i \left[(\lambda - w(i))\id - \Lambda^{>\ell}\right].
\end{equation}
Since $\left(\phi^{<\ell}\right)^\dagger L^{<\ell}\phi^{<\ell}$ is Hermitian, there exists a unitary $U_\ell$ and a diagonal matrix with non-decreasing diagonal entries $\Lambda^{<\ell}$, such that $\left(\phi^{<\ell}\right)^\dagger L^{<\ell}\phi^{<\ell} = U^\dagger_\ell \Lambda^{<\ell} U_\ell$.
Hence, we can write eq.~\eqref{eq:induction eval eq} as
\begin{equation}\label{eq:induction eval eq transformed}
    \Lambda^{<\ell} \left(\tilde C_\ell\right)_i = \left(\tilde C_\ell\right)_i \left[(\lambda - w(i))\id - \Lambda^{>\ell}\right] \quad \forall i\in[p],
\end{equation}
with $\tilde C_\ell = U_\ell C_\ell$.
Furthermore, for all $i\in[p]$
\begin{equation}\label{eq:inductive hypothesis proof}
    \begin{tikzpicture}[
        dot/.style = {circle,fill,minimum size=#1,inner sep=0pt,outer sep=0pt},
        dot/.default = 5pt,
        ]
        \node [dot,label=above:$\tilde C_\ell$] (cL) {};
        \node [dot,label=above:$\phi^{>\ell-1}$] (lhs) [right=.7 of cL] {};
        \node (mid) [right=.23 of cL] {};
        \node (topL) [below=.35 of cL] {};
        \node (topR) [below=.35 of lhs] {};
        \node (botL) [below=.03 of topL] {};
        \node (botR) [below=.03 of topR] {};
        \node [rectangle,draw] (Lleft) [below=.3 of mid] {$L^{>\ell-1}$};
        \node (left) [left=.3 of cL] {};
        \node (i) [below=.3 of botL] {$i$};
        \node (belowR) [below=.3 of botR] {};

        \draw (left) -- (cL) -- (lhs);
        \draw (cL) -- (topL);
        \draw (botL) -- (i);
        \draw[double] (lhs) -- (topR);
        \draw[double] (botR) -- (belowR);

        \node (eq) [right=.3 of lhs] {$=$};
        \node [dot,label=above:$\tilde C_\ell$] (cl) [right=.5 of eq] {};
        \node (left) [left=.3 of cl] {};
        \node (i) [below=.3 of cl] {$i$};
        \node [dot,label=above:$\phi^{>\ell}$] (phi) [right=.5 of cl] {};
        \node [dot,label=right:$L^{>\ell}$] (L) [below=.2 of phi] {};
        \node (below) [below=.2 of L] {};

        \draw (left) -- (cl) -- (phi);
        \draw[double] (phi) -- (L) -- (below);
        \draw (i) -- (cl);

        \node (plus) [right=.8 of phi] {$+$};
        \node [dot,label=above:$\tilde C_\ell$] (c) [right=.8 of plus] {};
        \node [dot,label=above:$\phi^{>\ell}$] (phi) [right=.5 of c] {};
        \node [dot,label=left:$\Omega$] (omega) [below=.2 of c] {};
        \node (i) [below=.2 of omega] {$i$};
        \node (below) [below=.5 of phi] {};
        \node (left) [left=.3 of c] {};
        
        \draw (i) -- (omega) -- (c) -- (phi);
        \draw[double] (phi) -- (below);
        \draw (left) -- (c);

        \node (eq) [below=1.5 of eq] {$=$};
        \node [dot,label=above:$\tilde C_\ell$] (cl) [right=.5 of eq] {};
        \node (left) [left=.3 of cl] {};
        \node (i) [below=.3 of cl] {$i$};
        \node [dot,label=above:$\Lambda^{>\ell}$] (lambda) [right=.5 of cl] {};
        \node [dot,label=above:$\phi^{>\ell}$] (phi) [right=.5 of lambda] {};
        \node (below) [below=.5 of phi] {};

        \draw (left) -- (cl) -- (phi);
        \draw[double] (phi) -- (below);
        \draw (i) -- (cl);

        \node (plus) [right=.8 of phi] {$+$};
        \node (omega) [right=.3 of plus] {$w(i)$};
        \node (left) [right=-.2 of omega] {};
        \node [dot,label=above:$\tilde C_\ell$] (c) [right=.3 of left] {};
        \node [dot,label=above:$\phi^{>\ell}$] (phi) [right=.5 of c] {};
        \node (i) [below=.3 of c] {$i$};
        \node (below) [below=.5 of phi] {};
        
        \draw (i) -- (c) -- (phi);
        \draw[double] (phi) -- (below);
        \draw (left) -- (c);

        \node (eq) [below=1.5 of eq] {$=$};
        \node (left) [right=.3 of eq] {};
        \node [dot,label=above:$\tilde C_\ell$] (c) [right=.3 of left] {};
        \node [dot,label=above:$\Lambda^{>\ell} + w(i)\id$] (middle) [right=1.5 of c] {};
        \node [dot,label=above:$\phi^{>\ell}$] (phi) [right=1.5 of middle] {};
        \node (i) [below=.3 of c] {$i$};
        \node (below) [below=.5 of phi] {};
        
        \draw (left) -- (c) -- (middle) -- (phi);
        \draw (i) -- (c);
        \draw[double] (phi) -- (below);
        
        \node (eq) [below=1.2 of eq] {$=$};
        \node (left) [right=.3 of eq] {};
        \node [dot,label=above:$\lambda\id - \Lambda^{<\ell}$] (op) [right=.8 of left] {};
        \node [dot,label=above:$\tilde C_\ell$] (c) [right=1.3 of op] {};
        \node [dot,label=above:$\phi^{>\ell}$] (phi) [right=.5 of c] {};
        \node (i) [below=.3 of c] {$i$};
        \node (below) [below=.3 of phi] {};
        \node (dot) [right=.5 of phi] {,};
        
        \draw (left) -- (op) -- (c) -- (phi);
        \draw (i) -- (c);
        \draw[double] (phi) -- (below);
    \end{tikzpicture}
\end{equation}
where in the last equality we have used eq.~\eqref{eq:induction eval eq transformed}.
After the gauge transformation $C_\ell \mapsto \tilde C_\ell, \ C_{\ell-1} \mapsto C_{\ell-1}U_\ell^\dagger$, we can write eq.~\eqref{eq:inductive hypothesis proof} in matrix notation as
\begin{equation}
    \phi^{>\ell-1} L^{>\ell-1} = \Lambda^{>\ell-1} \phi^{>\ell-1},
\end{equation}
where $\Lambda^{>\ell-1} = \lambda \id - \Lambda^{<\ell}$ has non-increasing entries.
Hence, we are left with a left canonical tensor train, such that the right interface vectors satisfy 
\begin{equation}
    \phi^{>\ell'}L^{>\ell'} = \Lambda^{>\ell'}\phi^{>\ell'} \quad \forall \ell'\ge\ell-1.
\end{equation}
This is the inductive hypothesis for $\ell -1$.

\emph{Conclusion.} We have found a gauge transformation that puts the tensor train into a form such that eq.~\eqref{eq:block sparsity 1 in proof} is satisfied.
We will now show that this in fact implies eq.~\eqref{eq:block sparsity 2 in proof}.
For any $\ell \in [d-1]$ we can write the eigenvalue equation as
\begin{equation}
    \begin{tikzpicture}[
        dot/.style = {circle,fill,minimum size=#1,inner sep=0pt,outer sep=0pt},
        dot/.default = 5pt
    ]
    \node (lmbd) {$\lambda$};
    \node [dot,label=above:$\phi^{<\ell+1}$] (c) [right=.3 of lmbd] {};
    \node (cb) [below=.5 of c] {};
    \node [dot,label=above:$\phi^{>\ell}$] (phi) [right=.7 of c] {};
    \node (phib) [below=.5 of phi] {};
    \node (eq) [right=.5 of phi] {$=$};

    \draw[double] (cb) -- (c) -- (phi);
    \draw[double] (phi) -- (phib);

    \node [dot,label=above:$\phi^{<\ell+1}$] (c2) [right=1.3 of eq] {};
    \node [dot,label=left:$L^{<\ell+1}$] (omega) [below=.2 of c2] {};
    \node (ob) [below=.2 of omega] {};
    \node [dot,label=above:$\phi^{>\ell}$] (phi2) [right=.7 of c2] {};
    \node (phi2b) [below=.5 of phi2] {};
    
    \draw[double] (ob) -- (omega) -- (c2);
    \draw (c2) -- (phi2);
    \draw[double] (phi2) -- (phi2b);
    
    \node (plus) [right=.5 of phi2] {$+$};
    \node [dot,label=above:$\phi^{<\ell+1}$] (c3) [right=.5 of plus] {};
    \node (c3b) [below=.5 of c3] {};
    \node [dot,label=above:$\phi^{>\ell}$] (phi3) [right=.7 of c3] {};
    \node [dot,label=right:$L^{>\ell}$] (L) [below=.2 of phi3] {};
    \node (Lb) [below=.2 of L] {};

    \draw[double] (c3b) -- (c3);
    \draw (c3) -- (phi3);
    \draw[double] (phi3) -- (L) -- (Lb);

    \node (eq2) [below=1.3 of eq] {$=$};
    \node [dot,label=above:$\phi^{<\ell+1}$] (c2) [right=1.3 of eq2] {};
    \node [dot,label=left:$L^{<\ell+1}$] (omega) [below=.2 of c2] {};
    \node (ob) [below=.2 of omega] {};
    \node [dot,label=above:$\phi^{>\ell}$] (phi2) [right=.7 of c2] {};
    \node (phi2b) [below=.5 of phi2] {};
    
    \draw[double] (ob) -- (omega) -- (c2);
    \draw (c2) -- (phi2);
    \draw[double] (phi2) -- (phi2b);

    \node (plus2) [right=.5 of phi2] {$+$};
    \node [dot,label=above:$\phi^{<\ell+1}$] (c4) [right=.5 of plus2] {};
    \node (c4b) [below=.5 of c4] {};
    \node [dot,label=above:$\Lambda^{>\ell}$] (Lmbd) [right=.7 of c4] {};
    \node [dot,label=above:$\phi^{>\ell}$] (phi4) [right=.5 of Lmbd] {};
    \node (phi4b) [below=.5 of phi4] {};
    \node (dot) [right=.3 of phi4] {,};

    \draw[double] (c4b) -- (c4);
    \draw (c4) -- (Lmbd) -- (phi4);
    \draw[double] (phi4) -- (phi4b);
    \end{tikzpicture}
\end{equation}
which, rearranging, can be written in matrix notation as 
\begin{equation}\label{eq:first core bs}
    L^{<\ell+1} \phi^{<\ell+1} \phi^{>\ell} = \phi^{<\ell+1} \left(\lambda\id - \Lambda^{>\ell}\right) \phi^{>\ell},
\end{equation}
which, since $\phi^{>\ell}$ has full row rank by the assumption of minimal ranks, implies
\begin{equation}
    L^{<\ell+1} \phi^{<\ell+1} = \phi^{<\ell+1} \left(\lambda\id - \Lambda^{>\ell}\right),
\end{equation}
which is eq.~\eqref{eq:block sparsity 2}.
$\square$

\section{Conditions on low-rankness}\label{app:conditions_on_lowrankness}

We here connect the conditions on low rank TT approximate representations of multivariate functions \cite{bachmayr_tensor_2016} with entanglement conditions on low rank TT approximations of quantum states \cite{schuch_entropy_2008,eisert_colloquium_2010}, known in this context as \emph{matrix product states} (MPS).

Multivariate $L^2(\mathbb{R}^d)$ functions can be naturally associated with $\ell^2(\mathbb{N}^d)$ sequences, where the $\ell^2$-norm is defined by $\|u\|_{\ell^2} \coloneqq \sqrt{\sum_{\mathbf{i}\in\mathbb{N}^d} u_\mathbf{i}^2}$.
Here, given a product basis $\Phi_{i_1,\dots,i_d}(x_1,\dots,x_d) = \Psi_{i_1}(x_1)\dots \Psi_{i_d}(x_d)$, where $\{\Psi_i\}_{i\in\mathbb{N}}$ is an orthonormal basis of $L^2(\mathbb{R})$, a function
\begin{equation}
    f(x_1,\dots,x_d) = \sum_{i_1,\dots,i_d\in\mathbb{N}} u_{i_1,\dots,i_d} \Phi_{i_1,\dots,i_d}(x_1,\dots,x_d) \quad \in L^2(\mathbb{R}^d)
\end{equation}
is associated with the sequence $(u_{\mathbf{i}})_{\mathbf{i}\in\mathbb{N}^d}$, in the sense that $\|f\|_{L^2} = \|u\|_{\ell^2}$. 

Let us introduce weak-$\ell^p$-norms of sequences, defined for $(a_n)_{n\in\mathbb{N}}$ by
\begin{equation}
    |a|_{w\ell^p} \coloneqq \sup_{n\in\mathbb{N}} n^{1/p} \tilde a_n,
\end{equation}
where $(\tilde a)_{n\in\mathbb{N}}$ is a non-increasing reordering of $(|a_n|)_{n\in\mathbb{N}}$.
These norms satisfy
\begin{equation}\label{eq:weak-lp-norm ineq}
    \|\cdot\|_{\ell^{p'}} \le |\cdot|_{w\ell^p} \le \|\cdot\|_{\ell^p}
\end{equation}
for $p'<p$.

Suppose a sequence $(u_\mathbf{i})_{\mathbf{i}\in\mathbb{N}^d} \in \ell^2(\mathbb{N}^d)$. 
For each $\eta \in [d-1]$ we define a linear operator $\mathcal{M}_\eta^u:\ell^2(\mathbb{N}^{d-\eta}) \rightarrow \ell^2(\mathbb{N}^\eta)$, which acts on $(v_\mathbf{i})_{\mathbf{i}\in\mathbb{N}^{d-\eta}}$ as
\begin{equation}
    \left(\mathcal{M}_\eta^u (v)\right)_{j_1,\dots,j_\eta} = \sum_{i_{1},\dots,i_{d-\eta}\in\mathbb{N}} u_{j_1,\dots,j_\eta,i_{1},\dots,i_{d-\eta}} v_{i_{1},\dots,i_{d-\eta}},
\end{equation}
which can be thought of as a matrification at the interface labeled $\eta$ of the tensor corresponding to $u$.
The operator $\mathcal{M}^u_\eta$ has a finite Frobenius norm $\|\mathcal{M}^u_\eta\|_F = \|u\|_{\ell^2}$, so its singular values $\sigma_\eta(u) \in \ell^2(\mathbb{N})$.
Let us now define the weak-$\ell^p_*$-norm for multi-index sequences, which we define for $u\in\ell^p(\mathbb{N}^d)$ (with $d\ge2$) by
\begin{equation}
    \|u\|_{w\ell^p_*} \coloneqq \max_{\eta\in[d-1]} |\sigma_\eta(u)|_{w\ell^p}.
\end{equation}

By \textit{Proposition 5.1} in ref.~\cite{bachmayr_tensor_2016}, the weak-$\ell^p_*$-norm determines the approximability of functions by low rank TTs.
Suppose $f\in L^2(\mathbb{R}^d)$, which is associated, through a given dictionary, with a sequence $u\in\ell^2(\mathbb{N}^d)$.
The Proposition states that, for $0<p<2$, there exists a sequence $\hat u \in \ell^2(\mathbb{N}^d)$, associated to a function $\hat f$, such that
\begin{equation} \label{eq:lowrankness_bound}
    \|f - \hat f\|_2 = \|u-\hat{u}\|_{\ell^2} \leq C \sqrt{d} \|u\|_{w\ell_*^p} \big(\max_{\eta \in [d-1]} \rank_\eta(\hat u)\big)^{-s},
\end{equation}
where $C>0$ is a constant, $s \coloneqq \frac1p - \frac12$ and $\rank_\eta(\hat u)$ counts the number of non-zero elements of $\sigma_\eta(\hat u)$.
To see the significance of this result, note that $\max_{\eta\in[d-1]}\rank_\eta(\hat u)$ is equal to the TT rank of $\hat u$, so through eq.~\eqref{eq:lowrankness_bound} the weak-$\ell^p_*$-norm puts limits on the existence of low rank TT approximations of $f$.
We will now show how this result relates to \emph{entanglement entropy scaling} used in quantum physics as a condition on low-rank MPS approximations of quantum states \cite{schuch_entropy_2008,eisert_colloquium_2010}. 

Quantum states are normalized elements of a complex Hilbert space $\mathcal{H}$ (or more accurately rays in $\mathcal{H}$) associated to a given system. 
Many-body systems have a tensor product structure, where $\mathcal{H} = \bigotimes_{k\in[d]} \mathcal{L}_k$, where $\mathcal{L}_k$ is the local Hilbert space of a single constituent sub-system and $d$ is the number of sub-systems.
Suppose $\mathcal{L}_k = \mathcal{L}$ for all $k\in[d]$ and choose an orthonormal basis $(\phi_i)_i$ for $\mathcal{L}$. 
This allows us to represent many-body quantum states by tensor networks, in an analogous way to how we use them to represent multivariate functions, although for quantum states we need to allow the tensor networks to be complex.

Given a quantum state $\psi$, the corresponding density matrix $\rho$ is defined as the rank-one projector onto $\psi$. 
For each $\eta \in [d-1]$, we define the reduced density matrix $\rho_\eta = \tr_\eta \rho$, where $\tr_\eta$ denotes the partial trace over all $\mathcal{L}_\ell$ with $\ell>\eta$.
The R\'enyi entropy for $\alpha\in[0,\infty]$ is defined by the analytic continuation of 
\begin{equation}\label{eq:renyi entropy}
    S_\alpha (\rho_\eta) \coloneqq \frac{\log\tr(\rho_\eta^\alpha)}{1-\alpha}.
\end{equation}

Suppose the singular value decomposition (SVD) $\mathcal{M}^\psi_\eta = U\Sigma V^\dagger$, where $U,V$ are isometries and $\Sigma = \diag(\sigma_\eta(\psi))$.
Hence, we can write $\rho_\eta = U\Sigma^2 U^\dagger$, which shows that the eigenvalues of $\rho_\eta$ are squares of the singular values $\sigma_\eta(\psi)$.
Therefore, we have
\begin{equation}\label{eq:entropy and lp}
    \tr(\rho_\eta^\alpha) = \|\sigma_\eta(\psi)\|_{\ell^{2\alpha}}^{2\alpha}.
\end{equation}

In quantum many-body physics, we are interested in approximations $(\hat \psi_d)_d$ of families of states $(\psi_d)_d$ on growing number $d$ of sub-systems, that satisfy $\|\psi_d - \hat \psi_d\|_2 \le \delta$ for all $d$.
We call such approximations efficient it the TT rank of $\hat \psi_d$ scales at most polynomially with $d$.
In \cite{schuch_entropy_2008} and \cite{verstraete_matrix_2006} it is shown that efficient approximations exist if, for some $0<\alpha<1$, there exist $c,c'>0$ such that for all $\eta\in[d-1]$ we have that $S_\alpha(\rho_{d,\eta}) = c\log d + c'$, where $\rho_{d,\eta} = \tr_\eta \rho_d$ and $\rho_d$ is the density matrix corresponding to $\psi_d$. 
We will now show that this follows from eq.~\eqref{eq:lowrankness_bound}.

Eq.~\eqref{eq:lowrankness_bound} implies that polynomial bond dimension approximations are guaranteed to exist if 
\begin{equation}
    \|\psi_d\|_{w\ell^p_*} \le \frac\delta{\sqrt{d}}(\poly(d))^s
\end{equation}
for some $0<p<2$.
Using eq.~\eqref{eq:weak-lp-norm ineq} and with $\alpha \coloneqq p/2$ we find a simpler sufficient condition 
\begin{equation}
    \max_{\eta\in[d-1]}\|\sigma_\eta(\psi_d)\|_{\ell^{2\alpha}} \le \frac\delta{\sqrt{d}}(\poly(d))^s.
\end{equation}
Through eq.~\eqref{eq:renyi entropy} and~\eqref{eq:entropy and lp}, we can rewrite this condition, for $0<\alpha<1$, as
\begin{equation}
    \max_{\eta\in[d-1]} S_\alpha(\rho_{\eta,d}) \le \frac{2\alpha}{1-\alpha} \log \frac\delta{\sqrt{d}}\left(\poly(d)\right)^s
\end{equation}
and hence there exist $c,c'>0$, such that $S_\alpha(\rho_{\eta,d}) \le c\log d + c'$ for all $\eta \in [d-1]$, which is the result of \cite{schuch_entropy_2008} and \cite{verstraete_matrix_2006}.

\end{document}